\documentclass{article}
\usepackage[T1]{fontenc}
\usepackage{amsmath}
\usepackage{amsthm}
\usepackage{amssymb}
\usepackage{hyperref}
\usepackage{listings}
\usepackage{float}
\usepackage{geometry}
\usepackage{tikz-cd}

\usepackage{graphicx,xcolor}

\newcommand{\bqn}{\begin{eqnarray}}
\newcommand{\eqn}{\end{eqnarray}}
\newcommand{\bq}{\begin{eqnarray*}}
\newcommand{\eq}{\end{eqnarray*}}

\newtheorem{example}{Example}[section]
\newtheorem{theorem}{Theorem}[section]

\begin{document}
\title{Statistical Inference in Persistent Homology (PH-STAT)}
%
%
\author{Moo K. Chung\\
Department of Biostatistics and Medical Informatics\\
University of Wisconsin-Madison\\
{\tt mkchung@wisc.edu}}

\date{Version 3, July 12, 2026}
\maketitle              
%


\begin{abstract}
We introduce PH-STAT, a comprehensive MATLAB toolbox designed for performing a wide range of statistical inference and learning tasks on persistent homology, primarily for graph and network data. Persistent homology is a prominent tool in topological data analysis (TDA) that captures the underlying topological features of complex data sets. The toolbox aims to provide users with an accessible and user-friendly interface for analyzing and interpreting topological data. The Matlab package is distributed in
\url{https://github.com/laplcebeltrami/PH-STAT}. This manual is downloaded from \url{https://arxiv.org/abs/2304.05912}. 
\end{abstract}

\tableofcontents

\section{Introduction}
\label{Intro}
PH-STAT (Statistical Inference on Persistent Homology) contains various statistical methods and algorithms for the analysis of persistent homology. The toolbox is designed to be compatible with a range of input data types, including point clouds, graphs, time series and functional data, graphs and networks, and simplicial complexes, allowing researchers from diverse fields to analyze their data using topological methods. The toolbox can be accessed and downloaded from the GitHub repository at \url{https://github.com/laplcebeltrami/PH-STAT}. The repository includes the source code, a detailed user manual, toy examples and example data sets for users to familiarize themselves with the toolbox's functionality. The user manual provides a comprehensive guide on how to run the toolbox, as well as an explanation of the underlying statistical methods and algorithms.

PH-STAT includes a rich set of statistical tools for analyzing and interpreting persistent homology, enabling users to gain valuable insights into their data. The toolbox provides visualization functions to help users understand and interpret the topological features of their data, as well as to create clear and informative plots for presentation and publication purposes. PH-STAT is an open-source project, encouraging users to contribute new features, algorithms, and improvements to the toolbox, fostering a collaborative and supportive community. PH-STAT is a versatile and powerful Matlab toolbox that facilitates the analysis of persistent homology in a user-friendly and accessible manner. By providing a comprehensive set of statistical tools, the toolbox enables researchers from various fields to harness the power of topological data analysis in their work.

\section{Morse Filtration}

In many applications, 1D functional data $f(t)$  is modeled as \cite{miller.1997}
\bqn f(t) = \mu(t) + \epsilon(t), \; t \in  \mathbb{R}, \label{eq:addictive} \eqn
where $\mu$ is the unknown mean signal to be estimated and $\epsilon$ is noise. In the usual statistical parametric mapping framework  \cite{friston.2002,kiebel.1996,worsley.1996}, inference on the model (\ref{eq:addictive}) proceeds as follows. If we denote an estimate of the signal by $\widehat \mu$, the residual $f - \widehat \mu$ gives an estimate of the  noise. One then constructs a test statistic $T(t)$, corresponding to a given hypothesis about the signal. As a way to account for temporal correlation of the statistic $T(t)$, the global maximum of the test statistic over the search space $\mathcal{M}$ is taken as the subsequent test statistic. Hence a great deal of the signal processing and statistical literature has been devoted to determining the distribution of $\sup_{t \in \mathbb{M}} T(t)$ using the random field theory \cite{taylor.2008,worsley.1996}, permutation tests \cite{nichols.2003} and the Hotelling--Weyl volume of tubes calculation
\cite {naiman.1990}. The use of the mean signal is one way of performing data reduction, however, this
may not necessarily be the best way to characterize complex multivariate imaging data.
Thus instead of using the mean signal, we can use persistent homology, which pairs local critical values \cite{edelsbrunner.2008,edelsbrunner.2002,zomorodian.2005}. It is intuitive that local critical values of $\widehat \mu$ approximately characterizes the shape of the continuous signal $\mu$ using only a finite number of scalar values. By pairing these local critical values in a nonlinear fashion and plotting them, one constructs the persistence diagram \cite{cohensteiner.2007,edelsbrunner.2008,morozov.2008,zomorodian.2001}.

\begin{figure}[t]
\begin{center}
\includegraphics[width=1\linewidth]{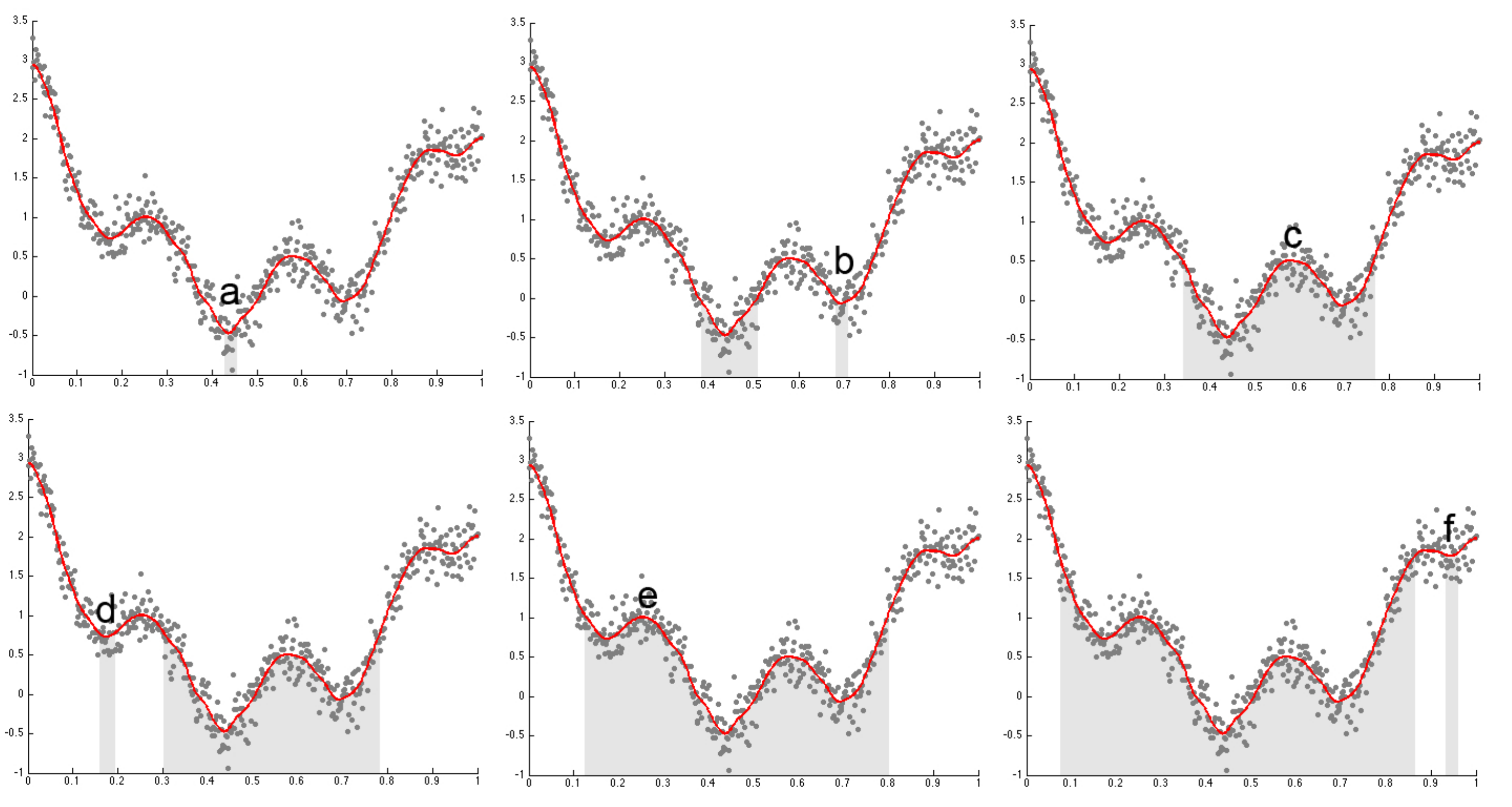}
\caption{\label{fig:bd} The births and deaths of connected components  in the sublevel sets in a Morse filtration \cite{chung.2009.IPMI}. We have local minimums $a < b < d <f$ and local maximums $c< e$. At $y=a$, we have a single connected component (gray area). As we increase the filtration value to $y=b$, we have the birth of a new component (second gray area). At the local maximum $y=c$, the two sublevel sets merge together to form a single component. This is viewed as the death of a component. The process continues till we exhaust all the critical values. Following the Elder rule, we pair birth to death: $(b,c)$ and $(d,e)$.  Other critical values are paired similarly. These paired points form the persistent diagram.}
\end{center}
\end{figure}

The function $\mu$ is called a {\em Morse function} if all critical values are distinct and non-degenerate, i.e., the Hessian does not vanish \cite{milnor.1973}. For a 1D Morse function $y=\mu(t)$, define sublevel set $R(y)$ as 
$$R_y =  \{t\in\mathbb{R}: \mu(t) \leq y\}.$$ 
As we increase height $y_1 \leq y_2 \leq  y_3 \leq\cdots$, the sublevel set gets bigger such that
$$R_{y_1} \subset R_{y_2} \subset R_{y_3} \subset \cdots.$$
The sequence of the sublevel sets form a {\em Morse filtration} with filtration values $y_1, y_2, y_3, \cdots$. 
Let $\beta_0(R_y)$ be the $0$-th Betti number of $R_y$, which counts the number of connected components in 
$R_y$. The number of connected components is the most often used topological invariant in applications \cite{edelsbrunner.2008}. $\beta_0(R_y)$ only changes its value as it passes through critical values (Figure \ref{fig:bd}). The birth and death of connected components in the Morse filtration is characterized by the pairing of local minimums and maximums. For 1D Morse functions, we do not have higher dimensional topological features beyond the connected components.

Let us denote the local minimums as $g_1, \cdots, g_m$ and the local maximums as $h_1, \cdots, h_n$. Since the critical values of a Morse function are all distinct, we can combine all minimums and maximums and reorder them from the smallest to the largest: 
We further order all critical values together and let
$$g_1 = z_{(1)} < z_{(2)} < \cdots < z_{(m+n)} = h_n,$$
where $z_{i}$ is either $h_i$ or $g_i$ and $z_{(i)}$ denotes the $i$-th largest number in $z_1, \cdots, z_{m+n}$.  In a Morse function, $g_1$ is smaller than $h_1$ and $g_m$ is smaller than $h_n$ in the unbounded domain $\mathbb{R}$ \cite{chung.2009.IPMI}. 

By keeping track of the birth and death of components, it is possible to compute topological invariants of sublevel sets such as the $0$-th Betti number $\beta_0$ \cite{edelsbrunner.2008}. As we move $y$ from $-\infty$ to $\infty$, at a local minimum, the sublevel set adds a new component so that 
$$\beta_0(R_{g_i - \epsilon}) = \beta_0( R_{g_i} ) +1$$
for sufficiently small $\epsilon$. This process is called the {\em birth} of the component. The newly born component is identified with the local minimum $g_i$. 

Similarly for at a local maximum, two components are merged as one so that
$$\beta_0(R_{h_i - \epsilon}) = \beta_0( R_{h_i} ) -1.$$
This process is called the {\em death} of the component. Since the number of connected components will only change if we pass through critical points and we can iteratively compute $\beta_0$ at each critical value as
$$\beta_0(R_{z_{(i+1)}}) = \beta_0(R_{z_{(i)}}) \pm 1.$$
The sign depends on if $z_{(i)}$ is maximum ($-1$) or minimum ($+1$). This is the basis of the Morse theory \cite{milnor.1973} that states that the topological characteristics of the sublevel set of Morse function are completely characterized by critical values.

To reduce the effect of low signal-to-noise ratio and to obtain smooth Morse function, either spatial or temporal smoothing have been often applied to brain imaging data before persistent homology is applied. In \cite{chung.2015.TMI,lee.2017.HBM}, Gaussian kernel smoothing was applied to 3D volumetric images. In \cite{wang.2018.annals}, diffusion was applied to temporally smooth data.

\begin{example}
As an example, elder's rule is illustrated in Figure \ref{fig:bd}, where the gray dots are simulated with Gaussian noise with mean 0 and variance $0.2^2$ as
$$f(x) = \mu(x) + N(0, 0.2^2)$$
with signal $\mu(t) = t + 7(t-1/2)^2 + \cos (7\pi t)/2$. The signal $\mu$ is estimated using heat kernel smoothing \cite{chung.2007.TMI} using degree $k=100$ and kernel bandwidth $\sigma=0.0001$ using {\tt WFS\_COS.m}. Now we  apply Morse filtration for filtration values $y$ from $-\infty$ to $\infty$. When we hit the first critical value $y=a$, the sublevel set consists of a single point. When we hit the minimum at $y=b$, we have the birth of a new component at $b$. When we hit the maximum at $y=c$, the two components identified by $a$ and $b$ are merged together to form a single component. When we pass through a maximum and merge two components, we pair the maximum with the higher of the two minimums of the two components \cite{edelsbrunner.2008}. Doing so we are pairing the birth of a component to its death. Obviously the paired extremes do not have to be adjacent to each other. If there is a boundary, the function value evaluated at the boundary is treated as a critical value. In the example, we need to pair $(b,c)$ and $(d,e)$. Other critical values are paired similarly. The  persistence diagram is then the scatter plot of these pairings computed using {\tt PH\_morse1D.m}. This is implemented as

\begin{verbatim}
t=[0:0.002:1]'; 
s= t + 7*(t - 0.5).^2 + cos(8*pi*t)/2; 
e=normrnd(0,0.2,length(x),1); 
Y=s+e; 

k=100; sigma=0.0001;
[wfs, beta]=WFS_COS(Y,x,k,sigma);

pairs=PH_morse1D(x,wfs);
\end{verbatim}
\end{example}

\section{Simplical Homology}
\index{simplical complex}

A high dimensional object can be approximated by the point cloud data $X$ consisting of $p$ number of points. If we connect points of which distance satisfy a given criterion, the connected points start to recover the topology of the object. Hence, we can represent the underlying topology as a collection of the subsets of $X$ that consists of nodes which are connected \cite{hart.1999,ghrist.2008,edelsbrunner.2010}. Given a point cloud data set $X$ with a rule for connections, the topological space is a simplicial complex and its element is a simplex \cite{zomorodian.2009}. For point cloud data, the Delaunay triangulation is probably the most widely used method for connecting points. The Delaunay triangulation represents the collection of points in space as a graph whose face consists of triangles. Another way of connecting point cloud data is based on Rips complex often studied in persistent homology. 

Homology is an algebraic formalism  to associate a sequence of objects with a topological space \cite{edelsbrunner.2010}.  In persistent homology, the algebraic formalism is usually built on top of objects that are hierarchically nested such as morse filtration, graph filtration and dendrograms. Formally homology usually refers to homology groups which are often built on top of a simplical complex for point cloud and network data \cite{lee.2014.MICCAI}.

The $k$-simplex $\sigma$ is the convex hull of $v+1$ independent points $v_0, \cdots, v_k$.  A point is a $0$-simplex, an edge is a $1$-simplex, and a filled-in triangle is a $2$-simplex. A {\em simplicial complex} is a finite collection of simplices such as points (0-simplex), lines (1-simplex), triangles (2-simplex) and higher dimensional counter parts \cite{edelsbrunner.2010}. A {\em $k$-skeleton} is a simplex complex  of up to $k$ simplices. Hence a graph is a 1-skeleton consisting of $0$-simplices (nodes) and $1$-simplices (edges). There are various simplicial complexes. The most often used simplicial complex in persistent homology is the Rips complex.

\begin{figure}[t]
\begin{center}
\includegraphics[width=1\linewidth]{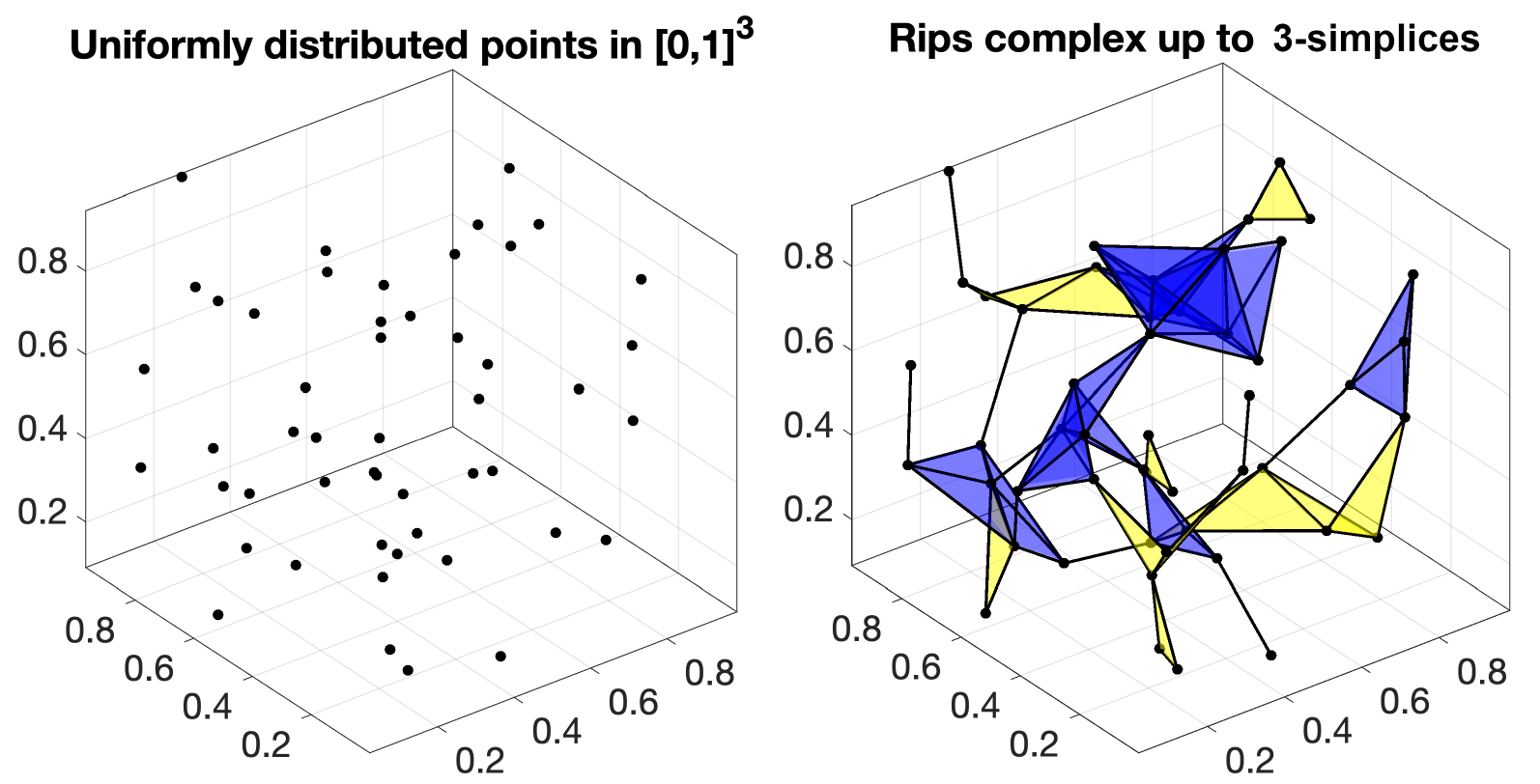}
\caption{Left: 50 randomly distributed points $X$ in $[0,1]^3$. Right: Rips complex $R_{0.3}(X)$ within radius 0.3 containing 106 1-simplices, 75 2-simplices (yellow) and 22 3-simplices (blue).}
\label{fig:rips}
\end{center}
\end{figure}

\subsection{Rips complex}

The Vietoris–Rips or Rips complex is the most often used simplicial complex in persistent homology. Let $X=\{ x_0, \cdots, x_p \}$ be the set of $n$ points in $\mathbb{R}^d$. The distance matrix between points in $X$ is given by $w=(w_{ij})$ where $w_{ij}$ is the distance between points $x_i$ and $x_j$. Then the Rips complex $R_\epsilon(X)$ is defined as follows \cite{edelsbrunner.2008,hatcher.2002}. The Rips complex is a collection of simplicial complexes parameterized by $\epsilon$. The complex $R_\epsilon(X)$ captures the topology of the point set $X$ at a scale of $\epsilon$ or less.

\begin{itemize}
\item The vertices of $R_\epsilon(X)$ are the points in $X$.

\item If the distance $w_{ij}$  is less than or equal to $\epsilon$, then there is an edge connecting points $x_i$ and $x_j$ in $R_\epsilon(X)$.

\item If the distance between any two points in $x_{i_0}, x_{i_1}, \ldots, x_{i_k}$ is less than or equal to $\epsilon$, then there is a $k$-simplex in $R_\epsilon(X)$ whose vertices are $x_{i_0}, x_{i_1}, \ldots, x_{i_k}$.

\end{itemize}

While a graph has at most 1-simplices, the Rips complex has at most $k$-simplices. In practice, the Rips complex is usually computed following the above definition  iteratively adding simplices of increasing dimension. Given $p+1$ number of points, there are potentially up to ${p+1 \choose k}$ $k$-simplices making the data representation extremely inefficient as the radius $\epsilon$ increases. Thus, we restrict simplices of dimension up to $k$ in practice. Such a simplical complex is called the $k$-skeleton. It is implemented as {\tt PH\_rips.m}, which inputs the matrix {\tt X} of size $p \times d$, dimension {\tt k} and radius {\tt e}. Then outputs the structured array {\tt S} containing the collection of nodes, edges,  faces up to $k$-simplices. For instance, the Rips complex up to 3-simplices in Figure \ref{fig:rips} is created using

\begin{verbatim}
p=50; d=3;
X = rand(p, d);
S= PH_rips(X, 3, 0.3)
PH_rips_display(X,S);

S =

  4×1 cell array

    { 50×1 double}      
    {106×2 double}      
    { 75×3 double}      
    { 22×4 double}      
\end{verbatim}

The Rips complex is then displayed using {\tt PH\_rips\_display.m} which inputs node coordinates {\tt X} and simplical complex {\tt S}.

\subsection{Boundary matrix}

Given a simplicial complex $K$, the boundary matrices $B_k$ represent the boundary operators between the simplices of dimension $k$ and $k-1$. Let $C_k$ be the collection of $k$-simplices. Define the  $k$-th boundary map $$\partial_k: C_k \rightarrow C_{k-1}$$ as a linear map that sends each $k$-simplex $\sigma$ to a linear combination of its $k-1$ faces
$$
\partial_k \sigma = \sum_{\tau \in F_k(\sigma)} (-1)^{\mathrm{sgn}(\tau, \sigma)} \tau,
$$
where $F_k(\sigma)$ is the set of $k-1$ faces of $\sigma$, and $\mathrm{sgn}(\tau, \sigma)$ is the sign of the permutation that sends the vertices of $\tau$ to the vertices of $\sigma$. This expression says that the boundary of a $k$-simplex $\sigma$ is the sum of all its $(k-1)$-dimensional faces, with appropriate signs determined by the orientation of the faces. The signs alternate between positive and negative depending on the relative orientation of the faces, as determined by the permutation that maps the vertices of one face to the vertices of the other face. The $k$-th boundary map  removes the filled-in interior of $k$-simplices. The vector spaces $C_k,$ $C_{k-1}, C_{k-2}, \cdots$ are then sequentially nested by boundary operator $\partial_k$ \cite{edelsbrunner.2010}:
 \bqn \cdots \xrightarrow{\partial_{k+1}} C_k \xrightarrow{\partial_k} C_{k-1} \xrightarrow{\partial_{k-1}} C_{k-2}  \xrightarrow{\partial_{k-2}}  \cdots. \label{eq:C_k} \eqn
 Such nested structure is called the {\em chain complex}.

\begin{figure}[t]
\centering
\includegraphics[width=0.8\linewidth]{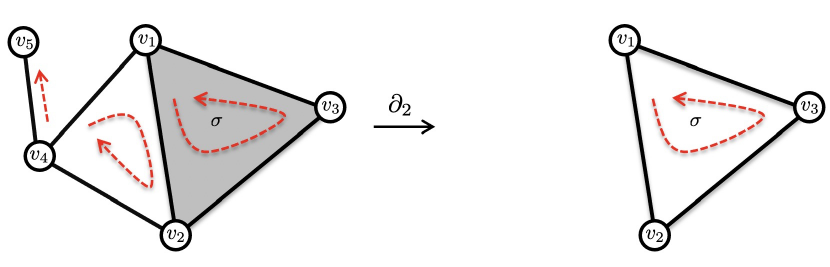}
\caption{A simplicial complex with 5 vertices and 2-simplex $\sigma =[v_1, v_2, v_3]$ with a filled-in face (colored gray). 
After boundary operation $\partial_2$, we are only left with 1-simplices $[v_1, v_2] + [v_2, v_3] + [v_3, v_1]$, which is the boundary of the filled in triangle. The complex has a single connected component ($\beta_0=1$) and a single 1-cycle. The dotted red arrows are the orientation of simplices.}
\label{fig:simplex}
\end{figure}

Consider a filled-in triangle $\sigma = [v_1, v_2, v_3] \in C_2$ with three vertices $v_1, v_2, v_3$ in Figure \ref{fig:simplex}. 
The boundary map $\partial_k$ applied to $\sigma$ resulted in the collection of three edges that forms the boundary of $\sigma$:
\bqn \partial_2 \sigma = [v_1,v_2] + [v_2,v_3] + [v_3,v_1] \in C_1. \label{eq: persist-p2} \eqn
If we give the direction or orientation to edges such that 
$$[v_3, v_1] = -[v_1, v_3],$$
and use edge notation $e_{ij} = [v_i, v_j]$, we can write (\ref{eq: persist-p2}) as
\bq \partial_2 \sigma = e_{12} + e_{23} + e_{31} = e_{12} + e_{23} - e_{13}. \eq

The boundary map can be represented as a boundary matrix $\boldsymbol{\partial}_k$ with respect to a basis of the vector spaces $C_k$ and $C_{k-1}$, where the rows of $\boldsymbol{\partial}_k$ correspond to the basis elements of $C_k$ and the columns correspond to the basis elements of $C_{k-1}$. The $(i,j)$ entry of $\boldsymbol{\partial}_k$ is given by the coefficient of the $j$th basis element in the linear combination of the $k-1$ faces of the $i$th basis element in $C_k$. The boundary matrix is  the higher dimensional version of the incidence matrix in graphs \cite{lee.2018.ISBI,lee.2019.CTIC,schaub.2018} showing how $(k-1)$-dimensional simplices are forming $k$-dimensional simplex. The $(i,j)$ entry of $\boldsymbol{\partial}_k$  is one if $\tau$ is a face of $\sigma$ otherwise zero. The entry can be -1 depending on the orientation of $\tau$. 
For the simplicial complex in Figure \ref{fig:simplex}, the boundary matrices are given by
\bq
\boldsymbol{\partial}_2 &=& 
\begin{array}{c}
\\
e_{12}\\
e_{23}\\
e_{31}\\
e_{24}\\
e_{41}\\
e_{45}
\end{array}
\begin{array}{c}
\sigma\\
\left(
\begin{array}{c}
 1 \\
 1   \\
 1 \\
 0\\
 0\\
 0 
\end{array}
\right)
\end{array}\\
\boldsymbol{\partial}_1 &=& 
\begin{array}{c}
\\
v_1\\
v_2\\
v_3\\
v_4\\
v_5
\end{array}
\begin{array}{c}
\begin{array}{cccccc}
e_{12} & e_{23} & e_{31} & e_{24} & e_{41} & e_{45}\\
\end{array}\\
\left(
\begin{array}{cccccc}
 -1 & 0 &  1  &  0&  1 & 0\\
 1 &  -1 & 0  &  -1&  0 &0\\
 0 &  1 & -1  &  0&  0 &0\\
 0 & 0& 0&  1&   -1 &-1\\
 0& 0& 0 & 0&  0 & 1
\end{array}
\right)
\end{array}\\
 \boldsymbol{\partial}_0 &=& 
\begin{array}{c}
\\
0
\end{array}
\begin{array}{c}
v_1 \;\; v_2 \;\; v_3 \;\;  v_4 \;\;  v_5\\
\left(
\begin{array}{ccccc}
 0 & 0 &  0 & 0 & 0
\end{array}
\right)
\end{array}.\\
\eq

\begin{figure}[t]
\begin{center}
\includegraphics[width=1\linewidth]{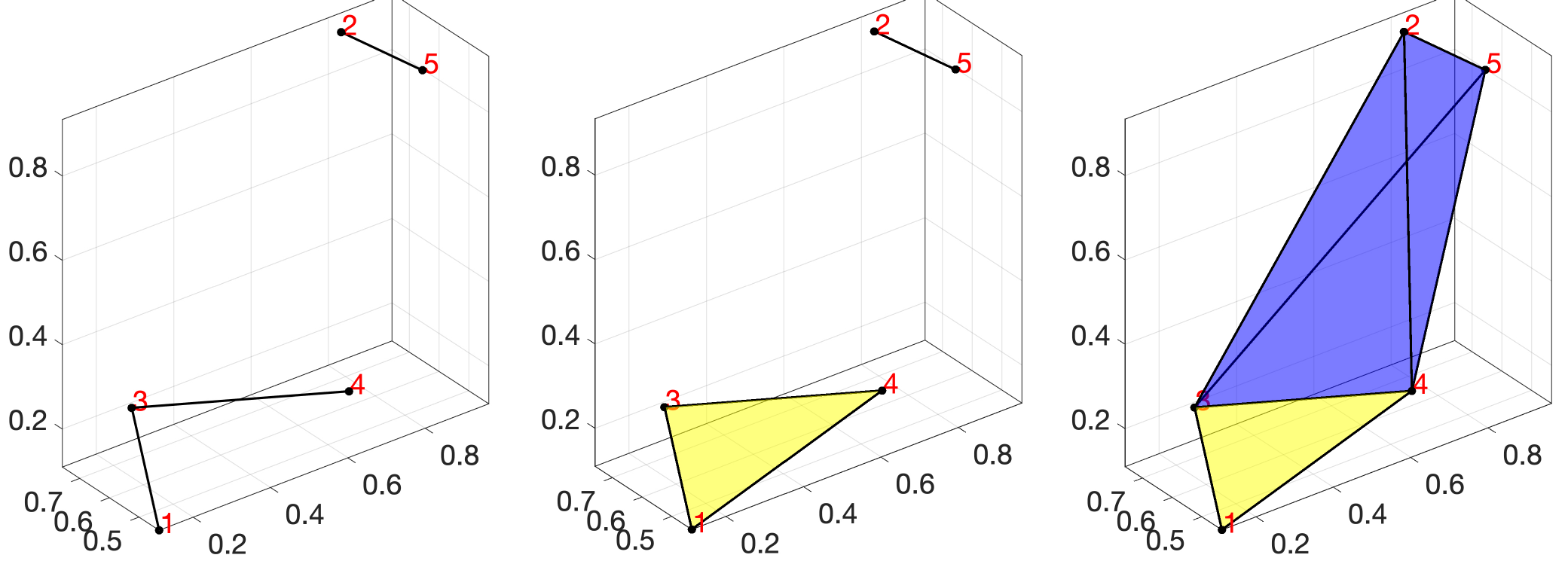}
\caption{Examples of boundary matrix computation. From the left to right, the radius is changed to 0.5, 0.6 and 1.0.}
\label{fig:boundary}
\end{center}
\end{figure}

In example in Figure \ref{fig:boundary}-left, {\tt PH\_rips(X,3,0.5)} gives
\begin{verbatim}
>> S{1}

     1
     2
     3
     4
     5
     
>> S{2}

     1     3
     2     5
     3     4
\end{verbatim}
{\tt PH\_boundary.m} only use node set {\tt S\{1\}} and edge set {\tt S\{2\}} in building boundary matrix {\tt B{1}} saving computer memory. 
\begin{verbatim}    
>> B{1}

    -1     0     0
     0    -1     0
     1     0    -1
     0     0     1
     0     1     0
\end{verbatim}
The columns of boundary matrix  {\tt B\{1\}} is indexed with the edge set in {\tt S\{2\}} such that the first column corresponds to edge {\tt [1,3]}. Any other potential edges {\tt [2,3]} that are not connected is simply ignored to save computer memory.

When we increase the filtration value and compute {\tt PH\_rips(X,3,0.6)}, a triangle is formed (yellow colored) and {\tt S\{3\}} is created (Figure \ref{fig:boundary}-middle). 
 
\begin{verbatim} 
>>  S{2}

     1     3
     1     4
     2     5
     3     4

>> S{3}

     1     3     4
\end{verbatim}
Correspondingly, the boundary matrices change to 
\begin{verbatim}
>> B{1}

    -1    -1     0     0
     0     0    -1     0
     1     0     0    -1
     0     1     0     1
     0     0     1     0

>> B{2}

     1
    -1
     0
     1
\end{verbatim}
From the edge set {\tt S\{2\}} that forms the row index for boundary matrix {\tt B\{2\}}, we have
{\tt [1,3] - [1,4] + [3,4]} that forms the triangle {\tt [1,3,4]}.

When we increase the filtration value further and compute {\tt PH\_rips(X,3,1)}, a tetrahedron is formed (blue colored) and {\tt S\{4\}} is created (Figure \ref{fig:boundary}-right). 
\begin{verbatim}
>> S{3}

     1     3     4
     2     3     4
     2     3     5
     2     4     5
     3     4     5

>>  S{4}

     2     3     4     5
\end{verbatim}     
 Correspondingly, the boundary matrix  {\tt B\{3\}} is created
\begin{verbatim}
>>  B{3}

     0
    -1
     1
    -1
     1
\end{verbatim}
The easiest way to check the computation is correct is looking at the sign of triangles in
{\tt - [2,3,4] + [2,3,5] - [2,4,5] + [3,4,5]}. Using the right hand thumb rule, which puts the orientation of triangle {\tt [3,4,5]} toward the center of the tetrahedron, the orientation of all the  triangles are toward the center of the tetrahedron. Thus, the signs are correctly assigned. Since computer algorithms are built inductively,  the method should work correctly in higher dimensional simplices. 
     
\subsection{Homology group}

The image of boundary map  is defined as
$$ \mathrm{im} \partial_{k+1} = \{ \partial_{k+1} \sigma | \sigma \in C_{k+1} \} \subset C_k,$$
which is a collection of boundaries. The elements of the image of $\partial_{k+1}$ are called $k$-boundaries, and they represent $k$-dimensional features that can be filled in by $(k+1)$-dimensional simplices. For instance, if we take the boundary $\partial_2$ of the triangle $\sigma$ in Figure \ref{fig:simplex}, we obtain a  1-cycle with edges $e_{12}, e_{23}, e_{31}$. The image of the boundary matrix  $B_{k+1}$ is the subspace spanned by its columns.  The column space can be found by the Gaussian elimination or singular value decomposition.

The kernel of boundary map is defined as
$$\mathrm{ker} \partial_{k} = \{ \sigma \in C_k | \partial_k \sigma = 0 \},$$
which is a collection of cycles. The elements of the kernel of $\partial_k$ are called cycles, since they form closed loops or cycles in the simplicial complex. The kernel of the boundary matrix $B_k$ is spanned by eigenvectors $v$ corresponding to zero eigenvalues of $B_k$. 

The boundary map satisfy the property that the composition of any two consecutive boundary maps is the zero map, i.e., 
\bqn \partial_{k-1} \circ \partial_k = 0. \label{eq:bb} \eqn
This reflect the fact that the boundary of a boundary is always empty. We can apply the boundary operation $\partial_1$ further to $\partial_2 \sigma$ in Figure \ref{fig:simplex} example and obtain 
\bq \partial_1 \partial_2 \sigma &=& 
\partial_1 e_{12}  + \partial_1 e_{23} + \partial_1 e_{31} \\
&=& v_2 - v_1 + v_3 - v_2 +  v_1 -v_3 =0.
\eq
This property (\ref{eq:bb}) implies that the image of $\partial_k$ is contained in the kernel of $\partial_{k-1}$, i.e., 
$$ \mathrm{im} \partial_{k+1} \subset \mathrm{ker} \partial_{k}.$$
Further, the boundaries $\mathrm{im} \partial_{k+1}$ form subgroups of the cycles $\mathrm{ker} \partial_{k}$. We can partition $\mathrm{ker} \partial_{k}$ into cycles that differ from each other by boundaries through the quotient space 
$$
H_k(K) = \mathrm{ker} \partial_k / \mathrm{im} \partial_{k+1},
$$
which is called the $k$-th homology group. $H_k(K)$  is a vector space that captures the $k$th topological feature or cycles in $K$. The elements of the $k$-th homology group are often referred to as $k$-dimensional cycles or $k$-cycles. Intuitively, it measures the presence of $k$-dimensional loops in the simplicial complex. 

\begin{figure}[t]
\begin{center}
\includegraphics[width=0.8\linewidth]{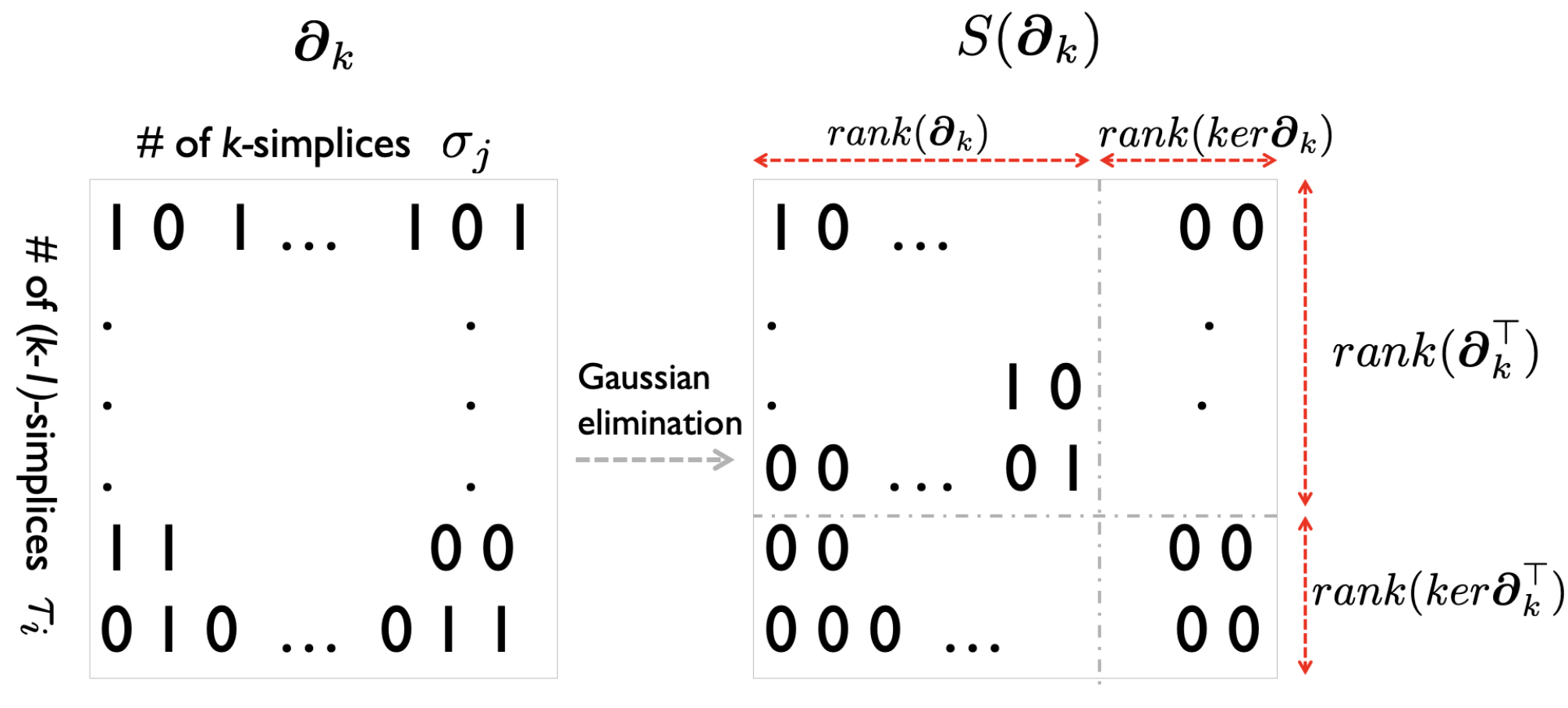}
\caption{The rank-nullity theorem for boundary matrix $\boldsymbol{\partial}_k$, which states the dimension of the domain of $\boldsymbol{\partial}_k$ is the sum of the dimension of its image and the dimension of its kernel (nullity). }
\label{fig:RN}
\end{center}
\end{figure}

The rank of $H_k(K)$ is the $k$th Betti number of $K$, which is an algebraic invariant that captures the topological features of the complex $K$. Although we put direction in the boundary matrices by adding sign, the Betti number computation will be invariant of how we orient simplices. The $k$-th Betti number $\beta_k$ is then computed as
\bqn \beta_k = rank (H_k) = rank (\mathrm{ker} \partial_k) - rank ( \mathrm{im} \partial_{k+1}). \label{eq:betti} \eqn 
The 0-th Betti number is the number of connected components while the 1-st Betti number is the number of cycles. The Betti numbers $\beta_k$ are usually algebraically computed by reducing the boundary matrix $\boldsymbol{\partial}_k$ to the Smith normal form $\mathcal{S}(\boldsymbol{\partial}_k)$, which has a block diagonal matrix as a submatrix in the upper left, via Gaussian elimination \cite{edelsbrunner.2010}. In the Smith normal form, we have the rank-nullity theorem for boundary matrix $\boldsymbol{\partial}_k$, which states the dimension of the domain of $\boldsymbol{\partial}_k$ is the sum of the dimension of its image and the dimension of its kernel (nullity) (Figure \ref{fig:RN}). In $\mathcal{S}(\boldsymbol{\partial}_k)$, the number of columns containing only zeros is $rank(\mathrm{ker} \boldsymbol{\partial}_k)$, the number of $k$-cycles while the number of columns containing one is $rank(\boldsymbol{\partial}_{k})$, the number of $k$-cycles that are boundaries.
Thus 
\bqn \beta_k = rank(ker \boldsymbol{\partial}_k) - rank (\boldsymbol{\partial}_{k+1}). \label{eq:rank}\eqn 
The computation starts with initial rank computation 
$$rank \boldsymbol{\partial}_0 = 0, \quad rank (ker \boldsymbol{\partial}_0) = p.$$
 
\begin{example} The boundary matrices  $\boldsymbol{\partial}_k$ in Figure \ref{fig:simplex} is transformed to the Smith normal form $\mathcal{S}(\boldsymbol{\partial}_k)$  after Gaussian elimination as
$$
\mathcal{S}(\boldsymbol{\partial}_1)= \left(
\begin{array}{cccccc}
 1 & 0 &  0  &  0&  0 & 0\\
 0 &  1 & 0  &  0&  0 &0\\
 0 &  0 & 1  &  0&  0 &0\\
 0 & 0& 0&  1&   0 & 0\\
 0& 0& 0 & 0&  0 & 0
\end{array}
\right), \quad  \mathcal{S}(\boldsymbol{\partial}_2)= \left(
\begin{array}{c}
 1 \\
 0   \\
 0 \\
 0\\
 0\\
 0 
\end{array}
\right).
$$
From (\ref{eq:rank}), the Betti number computation involves the rank computation of two boundary matrices. $rank(\boldsymbol{\partial}_0)=5$ is trivially the number of nodes in the simplicial complex. There are $rank(\mbox{ker}\boldsymbol{\partial}_1)=2$ zero columns and $rank (\boldsymbol{\partial}_1) = 4$ non-zero row columns. $rank(\boldsymbol{\partial}_2) =1$. Thus, we have
\bq \beta_0 &=& rank (\mbox{ker} \boldsymbol{\partial}_0) - rank (\boldsymbol{\partial}_1) =  5 - 4 = 1,\\
\beta_1 &=&  rank (\mbox{ker} \boldsymbol{\partial}_1) - rank (\boldsymbol{\partial}_2) =   2 - 1 = 1.
\eq 
\end{example}

\begin{figure}[t]
\begin{center}
\includegraphics[width=1\linewidth]{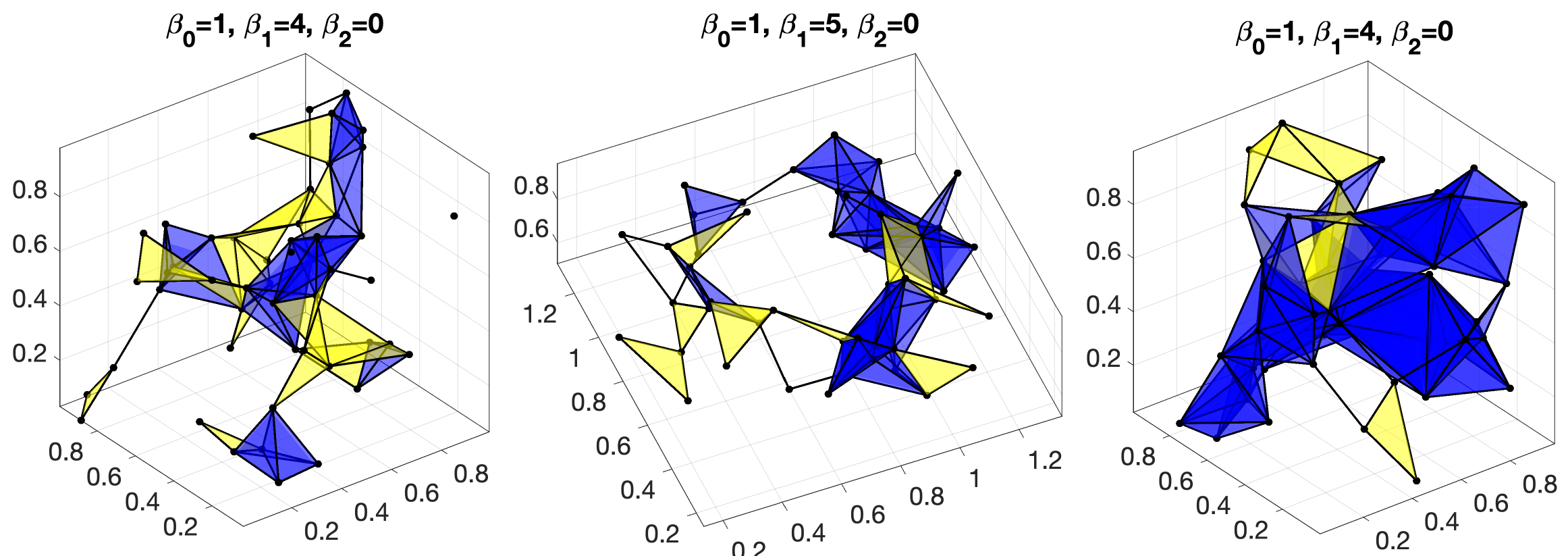}
\caption{Betti number computation on simplical complex using {\tt PH\_betti.m} function.}
\label{fig:betti}
\end{center}
\end{figure}

Following the above worked out example, the Betti number computation is implemented in {\tt PH\_boundary\_betti.m} which inputs the boundary matrices generated by {\tt PH\_boundary.m}. The function outputs $\beta_1, \beta_2, \cdots$. The function computes the $(d-1)$-th Betti number as 
\begin{verbatim}
betti(d)=  rank(null(B{d-1})) - rank(B{d}).
\end{verbatim}
Figure \ref{fig:betti} displays few examples of Betti number computation on Rips complexes. 
The rank computation in most computational packages is through the singular value decomposition (SVD).

\subsection{Rips filtrations} 
\index{Rips complex!triangulation}
\index{Rips}
\index{filtrations!Rips}
\index{filtration}

\begin{figure}[t]
\begin{center}
\includegraphics[width=0.7\linewidth]{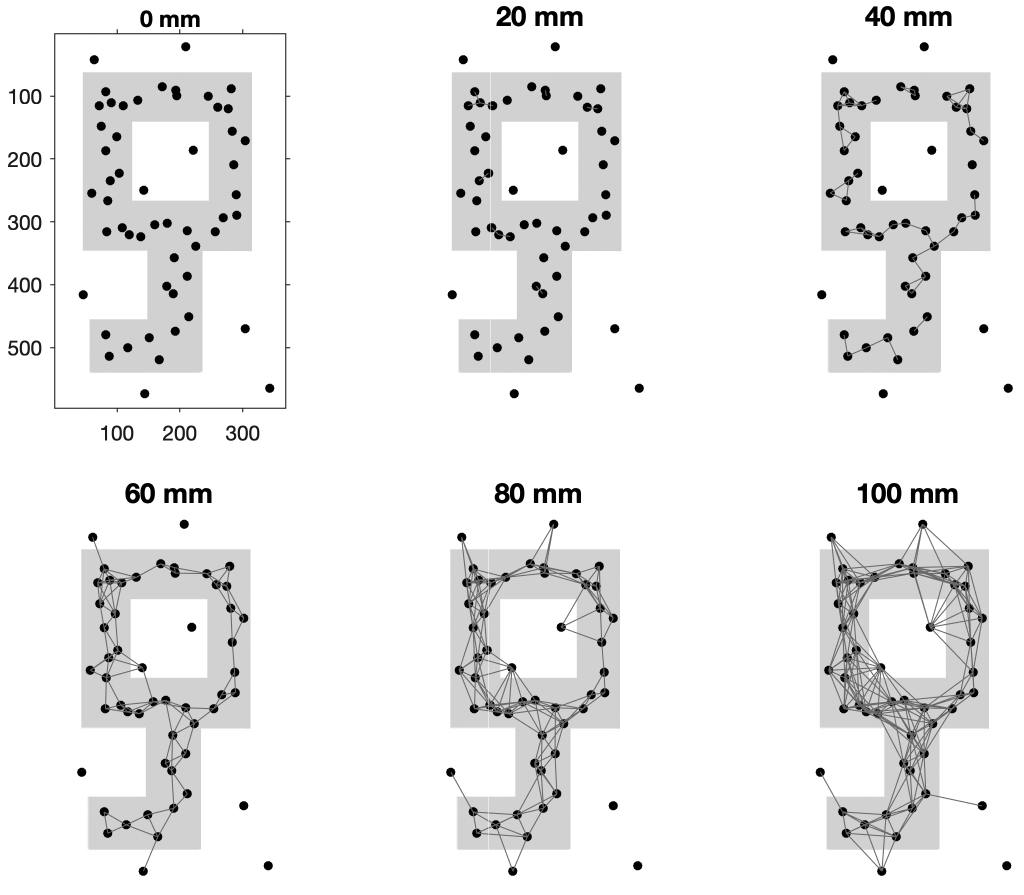}
\caption{Rips filtration on 1-skeleton of the point cloud data that was sampled along the underlying key shaped data. If two points are within the given radius, we connect them with an edge but do not form any other dimensional simplex. Such sparsity in Rips filtration can be more effective in practice. {\tt PH\_rips.m}  limit the dimension of skeleton we build Rips filtrations and do not build every possible simplical complexes.}
\label{fig:topology-ripskey}
\end{center}
\end{figure}

The Rips complex has the property that as the radius parameter value $\epsilon$ increases, the complex grows by adding new simplices. The simplices in the Rips complex at one radius parameter value are a subset of the simplices in the Rips complex at a larger radius parameter value. This nesting property is captured by the inclusion relation
$$\mathcal{R}{\epsilon_0} \subset \mathcal{R}{\epsilon_1}\subset \mathcal{R}{\epsilon_2}\subset \cdots $$
for $0=\epsilon_{0} \le \epsilon_{1} \le \epsilon_{2} \le \cdots.$
This nested sequence of Rips complexes is called the {\em Rips filtration}, which is the main object of interest in persistent homology (Figure \ref{fig:topology-ripskey}). The filtration values $\epsilon_0, \epsilon_1, \epsilon_2, \ldots$ represent the different scales at which we are studying the topological structure of the point cloud. By increasing the filtration value $\epsilon$, we are connecting more points, and therefore the size of the edge set, face set, and so on, increases.

The exponential growth in the number of simplices in the Rips complex as the number of vertices $p$ increases can quickly become a computational bottleneck when working with large point clouds. For a fixed dimension $k$, the number of $k$-simplices in the Rips complex grows as $\mathcal{O}(p^k)$, which can make computations and memory usage impractical for large values of $p$. Furthermore, as the filtration value $\epsilon$ increases, the Rips complex becomes increasingly dense, with edges between every pair of vertices and filled triangles between every triple of vertices. Even for moderately sized point clouds, the Rips filtration can become very ineffective as a representation of the underlying data at higher filtration values. The complex becomes too dense to provide meaningful insights into the underlying topological structure of the data. To address these issues, various methods have been proposed to sparsify the Rips complex. One such method is the graph filtration first proposed in \cite{lee.2011.ISBI,lee.2011.MICCAI}, which constructs a filtration based on a weighted graph representation of the data. The graph filtration can be more effective than the Rips filtration especially when  the topological features of interest are related to the graph structure of the data.


\subsection{Persistent diagrams}
 As the radius $\epsilon$ increases, the Rips complex  $R_{\epsilon}(X)$ grows and contains a higher-dimensional simplex that merges two lower-dimensional simplices representing the death of the two lower-dimensional features and the birth of a new higher-dimensional feature. The persistent diagram is a plot of the birth  and death  times of features. We start by computing the homology groups of each of the simplicial complexes in the filtration.

  \begin{figure}[t]
\center
  \includegraphics[width=1\textwidth]{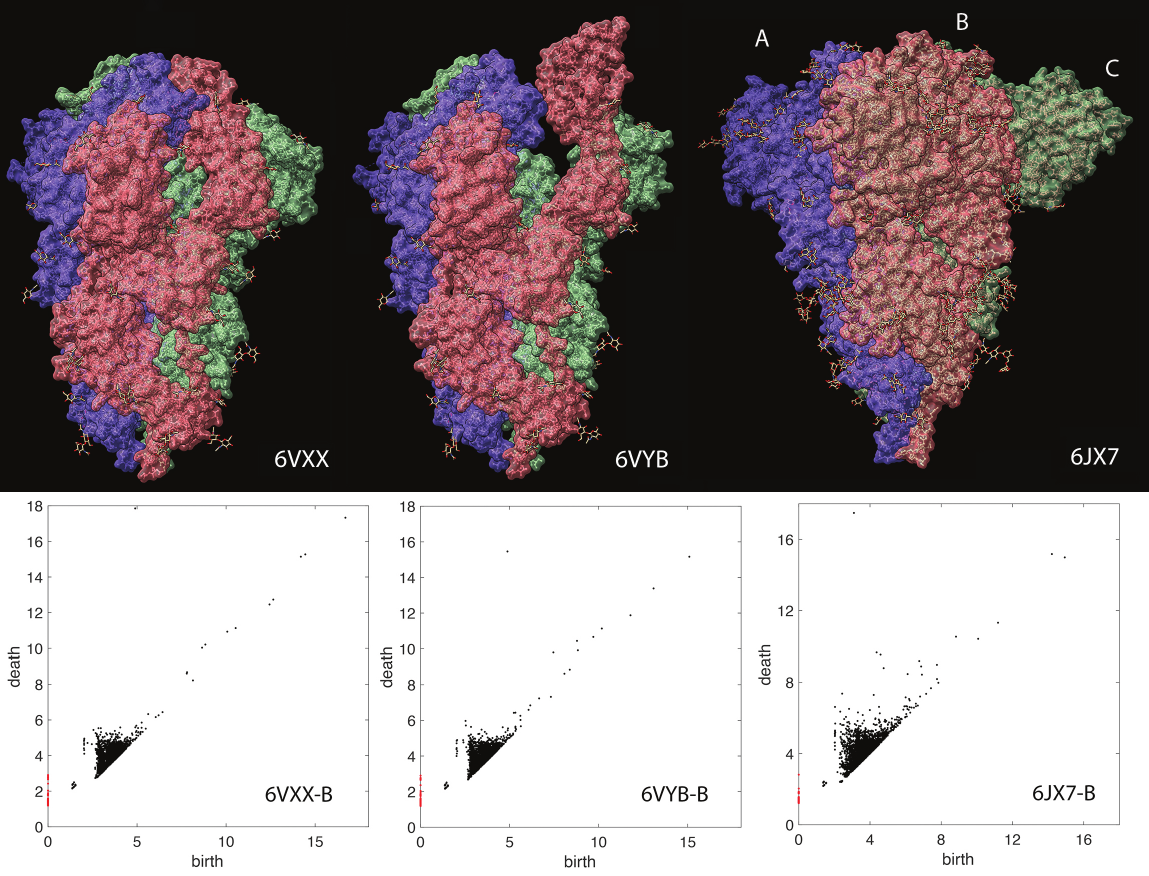}
  \caption{\small Top: Spike proteins of the three different corona viruses. The spike proteins consist of three similarly shaped interwinding substructures identified as A (blue), B (red) and C (green) domains. Bottom:  The persistent diagrams of spike proteins. The red dots are 0D homology and the black dots are 1D homology.}
  \label{fig:covid}
\end{figure}

 Let $H_k(K_i)$ denote the $k^{th}$ homology group of the simplicial complex $K_i$. We then track the appearance and disappearance of each homology class across the different simplicial complexes in the filtration. The birth time of a homology class  is defined as the smallest radius $\epsilon_b$ for which the class appears in the filtration, and the death time is the largest radius $\epsilon_d$ for which the class is present. We then plot each homology class as a point in the two-dimensional plane as $(\epsilon_b, \epsilon_d) \in \mathbb{R}^2$. The collection of all these points is the persistence diagram for $k$-th homology group.

To track the birth and death times of homology classes, we need to identify when a new homology class is born or an existing homology class dies as the radius $\epsilon$ increases. We can do this by tracking the changes in the ranks of the boundary matrices. Specifically, a $k$-dimensional cycle is born when it appears as a new element in the kernel of $\partial_k$ in a simplicial complex $K_i$ that did not have it before, and it dies when it becomes a boundary in $K_j$ for some $j > i$. Thus, we can compute the birth time $\epsilon_b$ of a $k$-dimensional homology class as the smallest radius for which it appears as a new element in the kernel of $\partial_k$. Similarly, we can compute the death time $\epsilon_d$ of the same class as the largest radius for which it is still a cycle in the simplicial complex $K_j$ for some $j>i$. By tracking the changes in the ranks of boundary matrices, we can compute the birth and death times of homology classes and plot them in the persistence diagram for the $k$-th homology group. However, the computation is fairly demanding and not scale well. 


\begin{example} The example came from \cite{chung.2021.MICCAI}.   The atomic structure of spike proteins of corona virus can be determined through the cryogenic electron microscopy (cryo-EM) \cite{cai.2020,walls.2020}. Figure \ref{fig:covid}-top displays a spike consists of three similarly shaped protein molecules with rotational symmetry often identified as A, B and C domains.  The 6VXX and 6VYB are respectively the closed and open states of SARS-Cov-2 from human \cite{walls.2020} while 6JX7 is feline coronavirus \cite{yang.2020}. We used the atomic distances in building Rips filtrations in computing persistent diagrams. The persistent diagrams of both closed and open states are almost identical in smaller birth and death values below 6 $\text{\normalfont\AA}$ (angstrom)  (Figure \ref{fig:covid}-bottom). The major difference is in the scatter points with larger brith and death values. However, we need a quantitative measures for comparing the topology of closed and open states. 
\end{example}

\section{Graph Filtration}

\subsection{Graph filtration}

The graph filtration has been the first type of filtrations applied in brain networks and it is now considered as the baseline filtrations in brain network data \cite{lee.2011.ISBI,lee.2011.MICCAI,lee.2012.TMI}. Euclidean distance is often used metric in building filtrations in persistent homology \cite{edelsbrunner.2010}. Most brain network studies  also use the Euclidean distances for building graph filtrations \cite{petri.2014,khalid.2014,cassidy.2015,wong.2016,anirudh.2016,palande.2017}. Given weighted network $\mathcal{X}=(V, w)$ with edge weight $w = (w_{ij})$, the binary network $\mathcal{X}_{\epsilon} =(V, w_{\epsilon})$ is a graph consisting of the node set $V$ and the binary edge weights 
$w_{\epsilon} =(w_{\epsilon,ij})$ given by 
\bqn w_{\epsilon,ij} =   \begin{cases}
1 &\; \mbox{  if } w_{ij} > \epsilon;\\
0 & \; \mbox{ otherwise}.
\end{cases}
\label{eq:case}
\eqn
Note \cite{lee.2011.MICCAI,lee.2012.TMI} defines the binary graphs by thresholding above such that $w_{\epsilon,ij} =1$ if $w_{ij} < \epsilon$ which is consistent with the definition of the Rips filtration. However, in brain connectivity,  higher value  $w_{ij}$ indicates stronger connectivity so we usually thresholds below  \cite{chung.2015.TMI}.

\begin{figure}[t]
\centering
\includegraphics[width=0.8\linewidth]{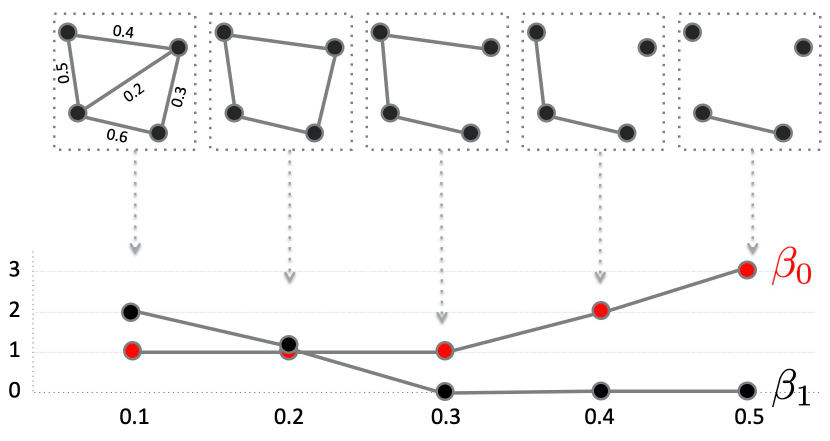}
\caption{Schematic of graph filtration and Betti curves. 
We sort the edge weights in an increasing order. We threshold the graph at filtration values and obtain binary graphs. The thresholding is performed sequentially by increasing the filtration values. The 0-th Betti number $\beta_0$, which counts the number of connected components, and the first Betti number $\beta_1$, which counts the number of cycles, is then plotted over the filtration. The Betti plots curves monotone in graph filtrations.}
\label{fig:persistent-betti01}
\end{figure}

\begin{figure}[t]
\centering
\includegraphics[width=1\linewidth]{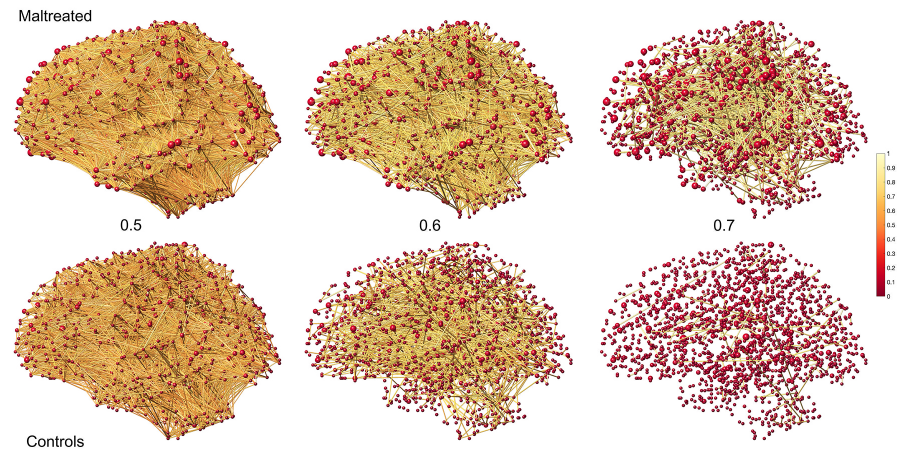}
\caption{Graph filtrations of maltreated children vs. normal control subjects on FA-values \cite{chung.2015.TMI}. The Pearson correlation is used as filtration values at 0.5, 0.6 and 0.7. maltreated subjects show much higher correlation of FA-values indicating more homogeneous and less varied structural covariate relationship.}
\label{fig:GF-maltreated}
\end{figure}

Note $w_{\epsilon}$ is the adjacency matrix of $\mathcal{X}_{\epsilon}$, which is a simplicial complex consisting of $0$-simplices (nodes) and $1$-simplices (edges)  \cite{ghrist.2008}. In the metric space $\mathcal{X}=(V, w)$, the Rips complex $\mathcal{R}_{\epsilon}(X)$ is a simplicial complex whose $(p-1)$-simplices correspond to unordered $p$-tuples of points that satisfy $w_{ij} \leq \epsilon$ in a pairwise fashion \cite{ghrist.2008}. While the binary network $\mathcal{X}_{\epsilon}$ has at most 1-simplices, the Rips complex can have at most $(p-1)$-simplices . Thus,  $\mathcal{X}_{\epsilon} \subset \mathcal{R}_{\epsilon}(\mathcal{X})$ and its compliment $\mathcal{X}_{\epsilon}^c \subset \mathcal{R}_{\epsilon}(\mathcal{X})$. Since a binary network is a special case of the Rips complex, we also have
$$\mathcal{X}_{\epsilon_0}  \supset \mathcal{X}_{\epsilon_1}  \supset \mathcal{X}_{\epsilon_2} \supset \cdots$$
and equivalently 
$$\mathcal{X}_{\epsilon_0}^c  \subset \mathcal{X}_{\epsilon_1}^c  \subset \mathcal{X}_{\epsilon_2}^c \subset \cdots  $$
for
$0=\epsilon_{0} \le \epsilon_{1} \le \epsilon_{2} \cdots.$ The sequence of such nested multiscale graphs  is defined as the {\em graph filtration} \cite{lee.2011.MICCAI,lee.2012.TMI}. Figure \ref{fig:persistent-betti01} illustrates a graph filtration in a 4-nodes example while Figure \ref{fig:GF-maltreated} shows the graph filtration on structural covariates on maltreated children on 116 parcellated brain regions. 

Note that $\mathcal{X}_0$ is the complete weighted graph while $\mathcal{X}_{\infty}$ is the node set $V$. By increasing the threshold value, we are thresholding at higher connectivity so more edges are removed. Given a weighted graph, there are infinitely many different filtrations. This makes the comparisons between two different graph filtrations difficult. 
For a graph with $p$ nodes, the maximum number of edges is $(p^2-p)/2$, which is obtained in a complete graph. If we order the edge weights in the increasing order, we have the sorted edge weights:
$$0 = w_{(0)} <  \min_{j,k} w_{jk} = w_{(1)} < w_{(2)} < \cdots < w_{(q)} = \max_{j,k} w_{jk},$$
where $q \leq (p^2-p)/2$.  The subscript $_{( \;)}$ denotes the order statistic. 
For all $\lambda < w_{(1)}$, $\mathcal{X}_{\lambda} = \mathcal{X}_{0}$ is the complete graph of $V$. For all $w_{(r)}  \leq \lambda < w_{(r+1)} \; (r =1, \cdots, q-1)$, $\mathcal{X}_{\lambda} = \mathcal{X}_{w_{(r)}}$. For all $ w_{(q)} \leq \lambda$, $\mathcal{X}_{\lambda}= \mathcal{X}_{\rho_{(q)}} =V$, the vertex set. Hence, the filtration given by
\bq  \mathcal{X}_{0}  \supset  \mathcal{X}_{w_{(1)}}  \supset  \mathcal{X}_{w_{(2)}}  \supset \cdots  \supset  \mathcal{X}_{w_{(q)}}\label{eq:maximal}\eq
is {\em maximal} in a sense that we cannot have any additional filtration $\mathcal{X}_{\epsilon}$ that is not one of the above filtrations. Thus, graph filtrations are usually given at edge weights \cite{chung.2015.TMI}.

The condition of having unique edge weights is not restrictive in practice. Assuming edge weights to follow some continuous distribution, the probability of any two edges being equal is zero. 
For discrete distribution, it may be possible to have identical edge weights. Then simply add Gaussian noise or add extremely small increasing sequence of numbers to $q$ number of edges.

\subsection{Monotone Betti curves}

The graph filtration can be quantified using monotonic function $f$ satisfying
\bqn f ( \mathcal{X}_{\epsilon_0} ) \geq f ( \mathcal{X}_{\epsilon_1} )  \geq f ( \mathcal{X}_{\epsilon_2} )  \geq \cdots   \label{eq:B} \eqn
or
\bqn f ( \mathcal{X}_{\epsilon_0} ) \leq f ( \mathcal{X}_{\epsilon_1} )  \leq f ( \mathcal{X}_{\epsilon_2} )  \leq \cdots   \label{eq:B2} \eqn

\begin{figure}[t]
\centering
\includegraphics[width=0.8\linewidth]{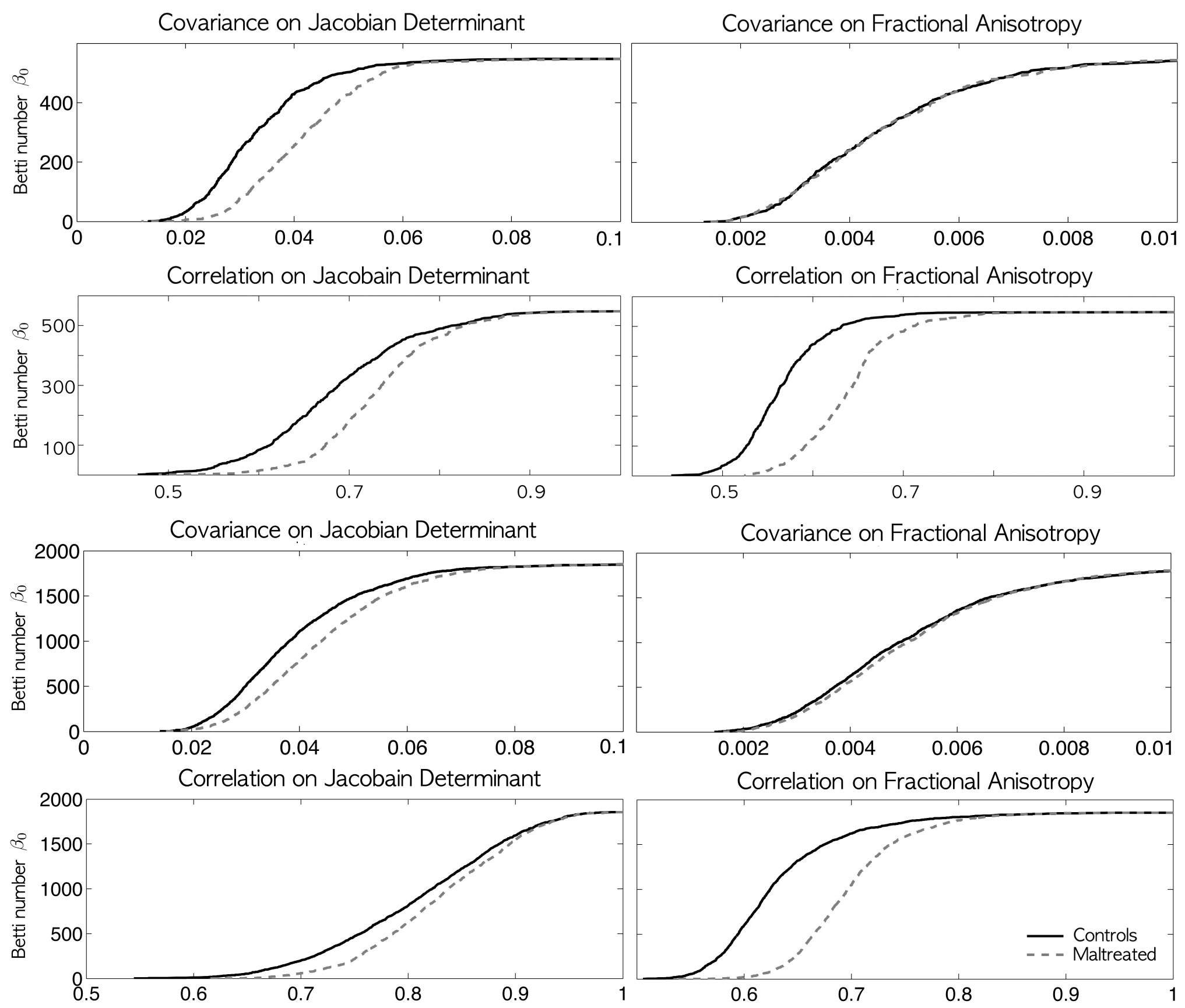}
\caption{The Betti curves on the covariance correlation matrices for Jacobian determinant (left column) and fractional anisotrophy (right column) on 548 (top two rows) and 1856 (bottom two rows) nodes \cite{chung.2015.TMI}. Unlike the covariance, the correlation seems to shows huge group separation between normal and maltreated children visually. However, in all 7 cases except top right (548 nodes covariance for FA), statistically significant differences were detected using the rank-sum test on the areas under the Betti-plots ($p$-value $< 0.001$).  The shapes of Betti-plots are consistent between the studies with different node sizes indicating the robustness of the proposed method over changing number of nodes.}
\label{fig:persist-betti0plots}
\end{figure}

The number of connected components (zeroth Betti number $\beta_0$) and the number of cycles (first Betti number $\beta_1$) satisfy the monotonicity (Figures \ref{fig:persistent-betti01} and \ref{fig:persist-betti0plots}). The size of the largest cluster also satisfies a similar but opposite relation of monotonic increase. There are numerous  monotone graph theory features \cite{chung.2015.TMI,chung.2017.IPMI}.

\begin{figure}[t]
\includegraphics[width=1\linewidth]{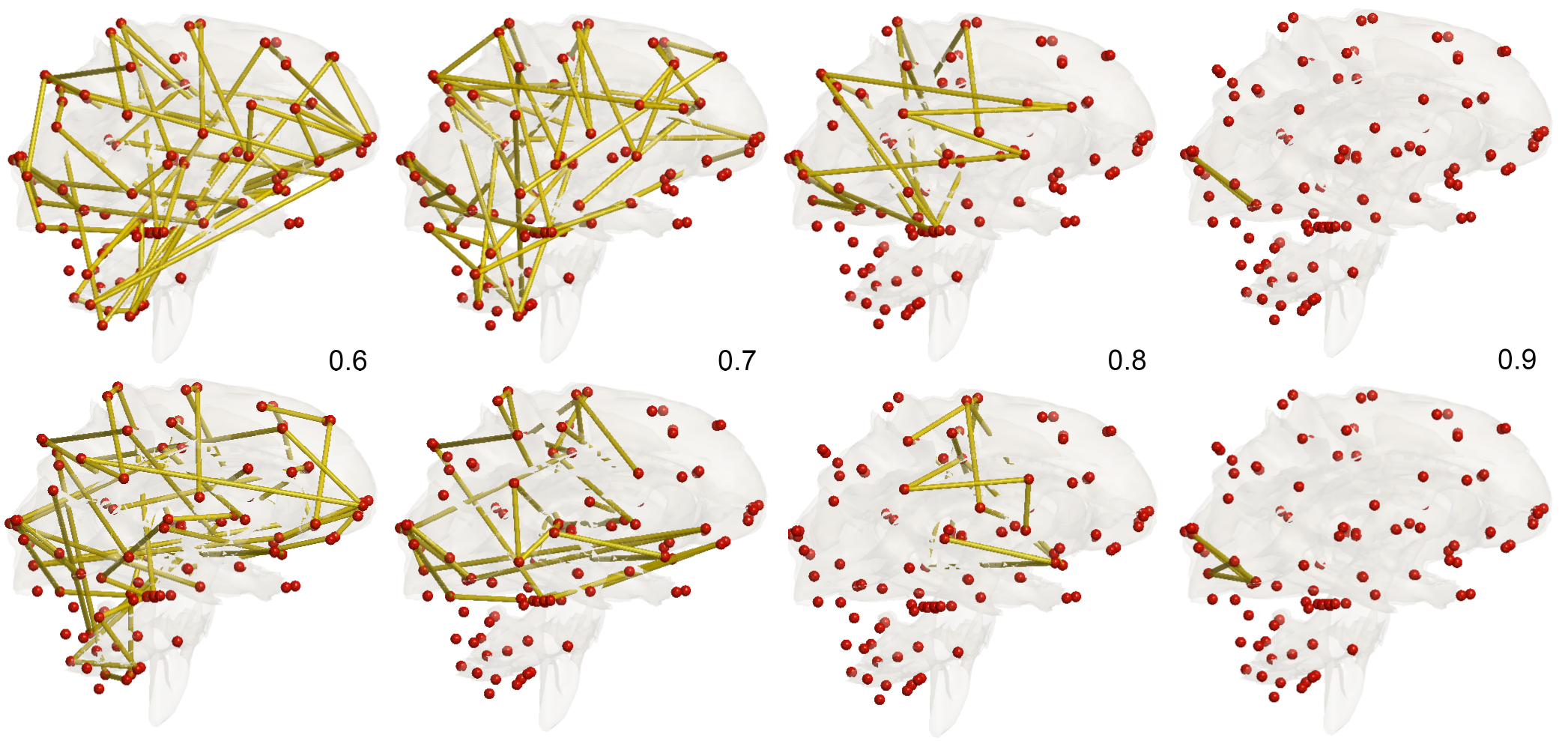}
\caption{The largest cycle at given correlation thresholds on rs-fMRI. Two representative subjects in HCP were used \cite{chung.2019.ISBI}. As the threshold increases, the length of cycles decreases monotonically.}
\label{fig:graphfiltration1}
\end{figure}

For graphs, $\beta_1$ can be computed easily as a function of $\beta_0$. Note that the Euler characteristic $\chi$ can be computed in two different ways
\bq \chi &=& \beta_0 - \beta_1 + \beta_2  - \cdots \\
 &=& \# nodes - \# edges + \# faces - \cdots,
 \eq
 where $\# nodes, \# edges, \# faces$ are the number of nodes, edges and faces. However, graphs do not have filled faces and Betti numbers higher than $\beta_0$ and $\beta_1$ can be ignored. Thus, a graph with $p$ nodes and $q$ edges is given by  \cite{adler.2010}
$$\chi = \beta_0 - \beta_1 = p - q.$$
Thus, $$\beta_1 = p - q - \beta_0.$$
In a graph, Betti numbers $\beta_0$ and $\beta_1$ are monotone over filtration on edge weights \cite{chung.2019.ISBI,chung.2019.NN}. When we do filtration on the maximal filtration in (\ref{eq:maximal}), edges are deleted one at a time.  
Since an edge has only two end points, the deletion of an edge disconnect the graph into at most two. Thus, the number of connected components ($\beta_0$) always increases and the increase is at most by one. Note $p$ is fixed over the filtration but $q$ is decreasing by one while $\beta_0$ increases at most by one. Hence, $\beta_1$ always decreases and the decrease is at most by one. Further, the length of the largest cycles, as measured by the number of nodes, also decreases monotonically (Figure \ref{fig:graphfiltration1}).

Identifying connected components in a network is important to understand in decomposing the network into disjoint subnetworks. The number of connected components (0-th Betti number) of a graph is a topological invariant that measures the number of structurally independent or disjoint subnetworks.  There are many available existing algorithms, which  are not related to persistent homology, for computing the number of connected components including the Dulmage-Mendelsohn decomposition \cite{pothen.1990}, which has been widely used for decomposing sparse matrices into block triangular forms in speeding up matrix operations.

In graph filtrations, the number of cycles increase or decreases as the filtration value increases. The pattern of monotone increase or decrease can visually show how the topology of the graph changes over filtration values. The overall pattern of  {\em Betti curves} can be used as a summary measure of quantifying how the graph changes over increasing edge weights \cite{chung.2013.MICCAI} (Figure \ref{fig:persistent-betti01}). The Betti curves are related to  barcodes. The Betti number is equal to the number of bars in the barcodes at the specific filtration value.

Figure \ref{fig:ripsvsgraph} displays an example of graph filtration 
constructed using random scatter points in a cube. Given scatter points {\tt X}, the pairwise distance matrix is computed as {\tt w = pdist2(X,X)}. The maximum distance is given by {\tt maxw = max (w(:))}. Betti curves are computed using {\tt PH\_betti.m}, which inputs the pairwise distance {\tt w} and the range of filtration values{\tt [0:0.05:maxw]}. The function outputs $\beta_0$ and $\beta_1$ values as structured arrays {\tt beta.zero} and {\tt beta.one}. We display them using {\tt PH\_betti\_display.m}.

\begin{verbatim}
p=50; d=3;
X = rand(p, d);
w = pdist2(X,X); 
maxw = max(w(:)); 

thresholds=[0:0.05:maxw];
beta = PH_betti(w, thresholds);
PH_betti_display(beta,thresholds)
\end{verbatim}

\subsection{Rips filtration vs. graph filtration}

\begin{figure}[t]
\begin{center}
\includegraphics[width=0.8\linewidth]{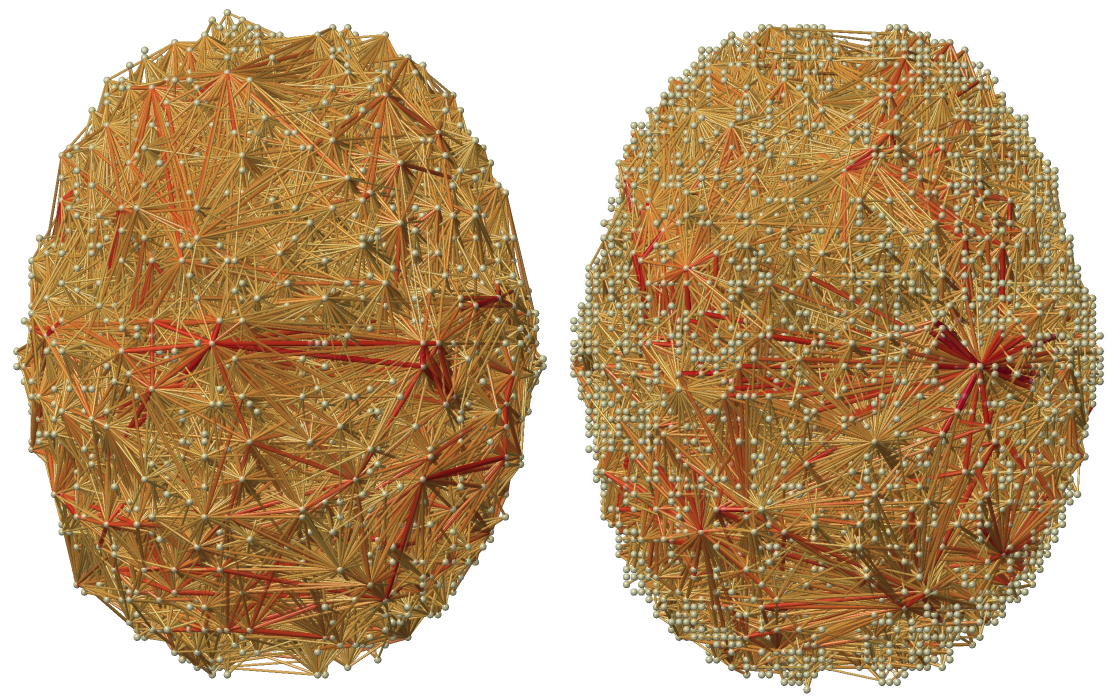}
\caption{rs-fMRI correlation network of two subjects from HCP with more than 25000 nodes. Identifying cycles and computing the number of cycles can be computationally demanding in this type of dense correlation network since persistent homology computations are not very scalable. }
\label{fig:persist-cycle}
\end{center}
\end{figure}

Persistent homology does not scale well with increased data size (Figure \ref{fig:persist-cycle}). The computational complexity of persistent homology grows rapidly with the number of simplices  \cite{topaz.2015}. With $p$ number of nodes, the size of the $k$-skeleton  grows as $p^{k+1}$. Homology calculations are often done by Gaussian elimination, and if there are $N$ simplices, it takes ${\cal O} (N^3)$ time to perform. In $\mathbb{R}^d$, the computational complexity is ${\cal O}(p^{3k+3})$ \cite{solo.2018}. Thus, the computation of Rips complex is exponentially costly. It can easily becomes infeasible when when one tries to use brain networks at the voxel level resolution. Thus, there have been many attempts in computing Rips complex approximately but fast for large-scale data such as
alpha filtration based on Delaunay triangulation with ${\cal O} (n^2)$ for $k=3$ and sparse Rips filtration with  ${\cal O}(n)$ simplicies and ${\cal O} (n \log n)$ runtime \cite{chung.2020.arXiv,otter.2017,sheehy.2013}. However, all of these filtrations are all approximation to Rips filtration. To remedy the computational bottleneck caused by Rips filtrations, {\em graph filtration} was introduced  particularly for network data \cite{lee.2011.MICCAI,lee.2012.TMI}.

\begin{figure}[t]
\includegraphics[width=1\linewidth]{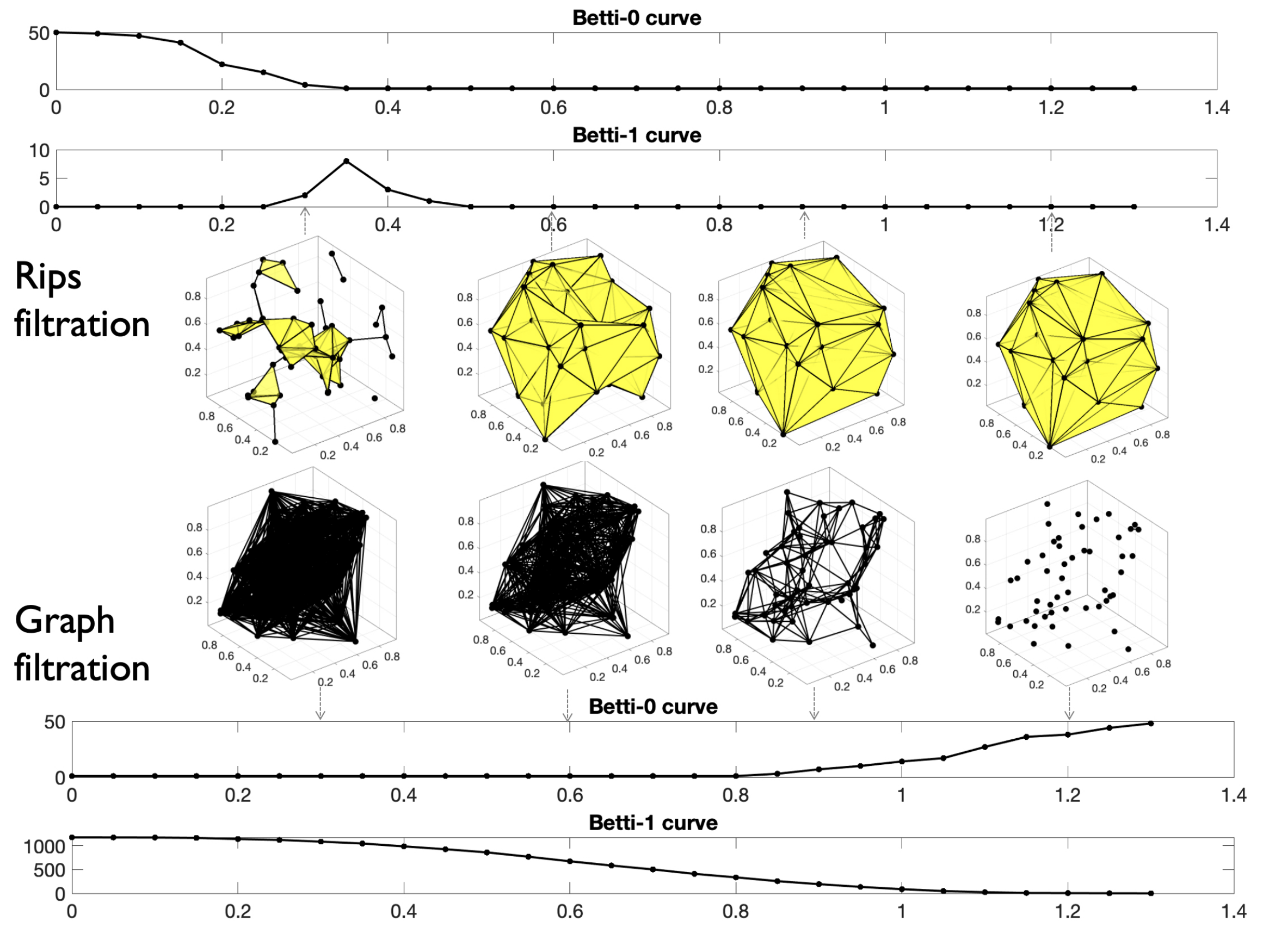}
\caption{The comparison between the Rips and graph filtrations.}
\label{fig:ripsvsgraph}
\end{figure}

The graph filtration is a special case of Rips filtration restricted to 1-simplices. If the Rips filtration up to 2-simplices is given by {\tt PH\_rips(X,2, e)}, the graph filtration is given by {\tt PH\_rips(X,1, maxw-e)}. In Figure \ref{fig:ripsvsgraph} displays the comparison of two filtrations for randomly generated 50 nodes in a cube. In the both filtrations, $\beta_0$-curves are monotone. However, $\beta_1$-curve for the Rips filtration is not monotone. Further, the range of changes in $\beta_1$ is very narrow. In some randomly generated points, we can have multiple peaks in $\beta_1$ making the $\beta_1$-curve somewhat unstable. On the other hand, the $\beta_1$-curve for the graph filtration is monotone and gradually changing over the whole range of filtration values. This will give consistent to the $\beta_1$-curve that is required for increasing statistical power in the group level inference.

\subsection{Graph filtration in trees}

Binary trees have been a popular data structure to analyze using persistent homology in recent years \cite{bendich.2016,li.2017}. Trees and graphs are 1-skeletons, which are Rips complexes consisting of only nodes and edges. However, trees do not have 1-cycles and can be quantified using up to 0-cycles only, i.e., connected components, and higher order topological features can be simply ignored. However, \cite{garside.2020} used somewhat inefficient filtrations in the 2D plane that increase the radius of circles from the root node or points along the circles. Such filtrations will produces persistent diagrams (PD) that spread points in 2D plane. Further, it may create 1-cycles. Such PD are difficult to analyze since scatter points do not correspond across different PD. For 1-skeleton, the {\em graph filtration} offers more efficient alternative \cite{chung.2019.NN,song.2023}.

\begin{figure}[t]
\centering
\includegraphics[width=1\linewidth]{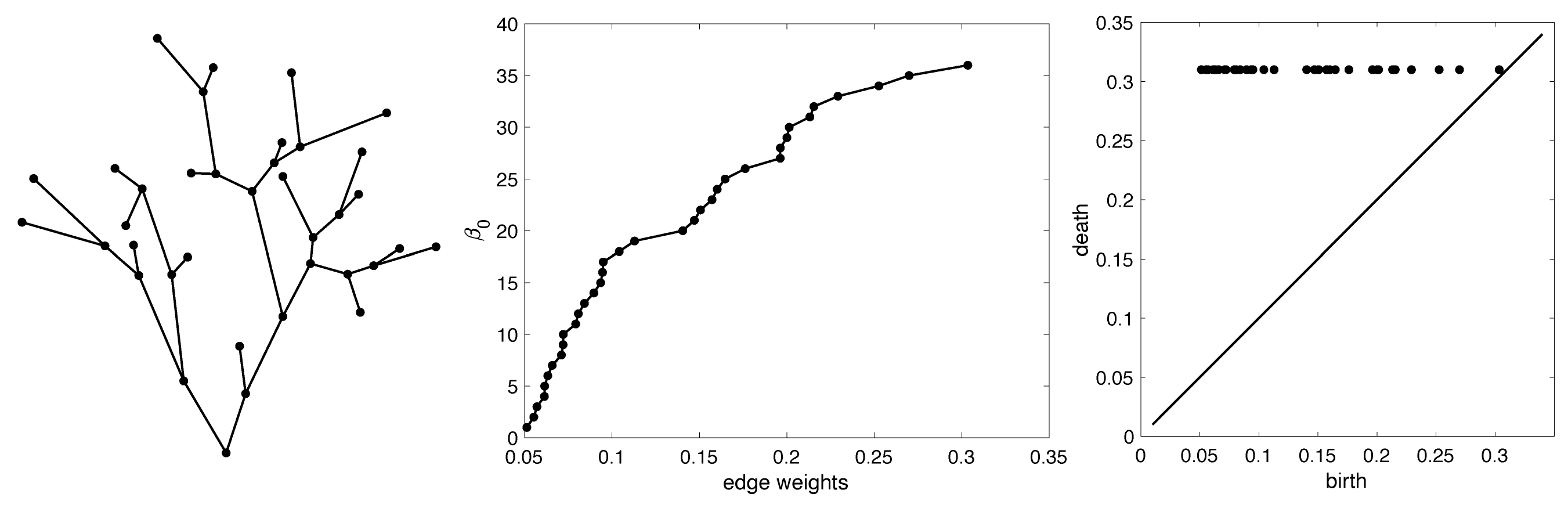}
\caption{Left: binary tree used in  \cite{garside.2020}. Middle: $\beta_0$-curve over graph filtration. Edge weights of the tree is used as the filtration values.  Right: The points in the persistent diagram all lined up at $y=0.31$, which is arbitarily picked to be larger than the maximum edge weight 0.3034.}
\label{fig:filtrationtree}
\end{figure}

Consider a tree $\mathcal{T}=(V, w)$ with node set $V=\{1, 2, \cdots, p \}$ and weighted adjacency matrix $w$. If we have binary tree with binary adjacency matrix, we add edge weights by taking the distance between nodes $i$ and $j$ as the edge weights $w_{ij}$ and build a weighted tree with $w=(w_{ij})$. For a tree $T$ with $p \geq 2$ nodes and unique $p-1$ positive edge weights $w_{(1)} < w_{(2)} < \cdots < w_{(p-1)}$. Threshold $\mathcal{T}$ at $\epsilon$ and define the binary tree $\mathcal{T}_{\epsilon} =(V, w_{\epsilon})$ with edge weights  $w_{\epsilon} = (w_{\epsilon,ij}), w_{\epsilon,ij} =  1$ if $w_{ij} > \epsilon$ and 0 otherwise. Then we have graph filtration on trees
\bqn  \mathcal{T}_{w_{(0)}}  \supset  \mathcal{T}_{w_{(1)}}  \supset  \mathcal{T}_{w_{(2)}}  \supset \cdots  \supset  \mathcal{T}_{w_{(p-1)}}.\label{eq:maximal2}\eqn
Since all the edge weights are above filtration value $w_{(0)}=0$, all the nodes are connected, i.e., 
$\beta_0(w_{(0)}) = 1$. Since no edge weight is above the threshold  $w_{(q-1)}$, $\beta_0( w_{(p-1)}) = p$. 
Each time we threshold, the tree splits into two and  the number of components increases exactly by one in the tree  \cite{chung.2019.NN}. Thus,
we have 
$$\beta_0(\mathcal{T}_{w_{(1)}}) = 2, \beta_0(\mathcal{T}_{w_{(2)}}) = 3, \cdots, \beta_0(\mathcal{T}_{w_{(p-1)}}) = p.$$ 
Thus, the coordinates for the 0-th Betti curve is given by
$$(0, 1), (w_{(1)}, 2), \cdots,  (w_{(2)}, 3), (w_{(p-1)}, p), (\infty, p).$$

All the 0-cycles (connected components) never die once they are bone over graph filtration. For convenience, we simply let the death value of  0-cycles  at some fixed number $c > w_{(q-1)}$. Then PD of the graph filtration is simply 
$$(w_{(1)},c), (w_{(2)}, c), \cdots, (w_{(q-1)}, c)$$ 
forming 1D scatter points along the horizontal line $y=c$ making various operations and analysis on PD significantly simplified  \cite{song.2023}. Figure \ref{fig:filtrationtree} illustrates the graph filtration and corresponding 1D scatter points in PD on the binary tree used in  \cite{garside.2020}. A different graph filtration is also possible by making the edge weight to be the shortest distance from the root node. However, they should carry the identical topological information.

For a general graph, it is not possible to analytically determine the coordinates for its Betti curves. The best we can do is to compute the number of connected components $\beta_0$ numerically using the single linkage dendrogram method (SLD) \cite{lee.2012.TMI}, the Dulmage-Mendelsohn decomposition \cite{pothen.1990,chung.2011.SPIE} or through the Gaussian elimination \cite{deSilva.2007,carlsson.2008,edelsbrunner.2002}.

\subsection{Birth-death decomposition}

Unlike the Rips complex, there are no higher dimensional topological features beyond the 0D and 1D topology in graph filtration. The 0D and 1D persistent diagrams  $(b_i, d_i)$ tabulate the life-time of 0-cycles (connected components) and 1-cycles (loops) that are born at the filtration value $b_i$ and die at value $d_i$. The 0th Betti number $\beta_0(w_{(i)})$ counts the number of 0-cycles at filtration value $w_{(i)}$  and shown to be non-decreasing over filtration (Figure \ref{fig:BDschematic}) \cite{chung.2019.ISBI}: 
$\beta_0(w_{(i)}) \leq \beta_0(w_{(i+1)}).$
On the other hand the 1st Betti number $\beta_1(w_{(i)})$ counts the number of independent loops and shown to be non-increasing over filtration  \cite{chung.2019.ISBI}:
$\beta_1(w_{(i)}) \geq \beta_1(w_{(i+1)}).$

During the graph filtration, when new components is born, they never dies. Thus, 0D persistent diagrams are completely characterized by birth values $b_i$ only. Loops are viewed as already born at $-\infty$. Thus, 1D persistent diagrams are completely characterized by death values $d_i$ only. We can show that the edge weight set $W$ can be partitioned into 0D birth values and 1D death values \cite{song.2021.MICCAI}:

\begin{figure}[t]
\begin{center}
\includegraphics[width=1\linewidth]{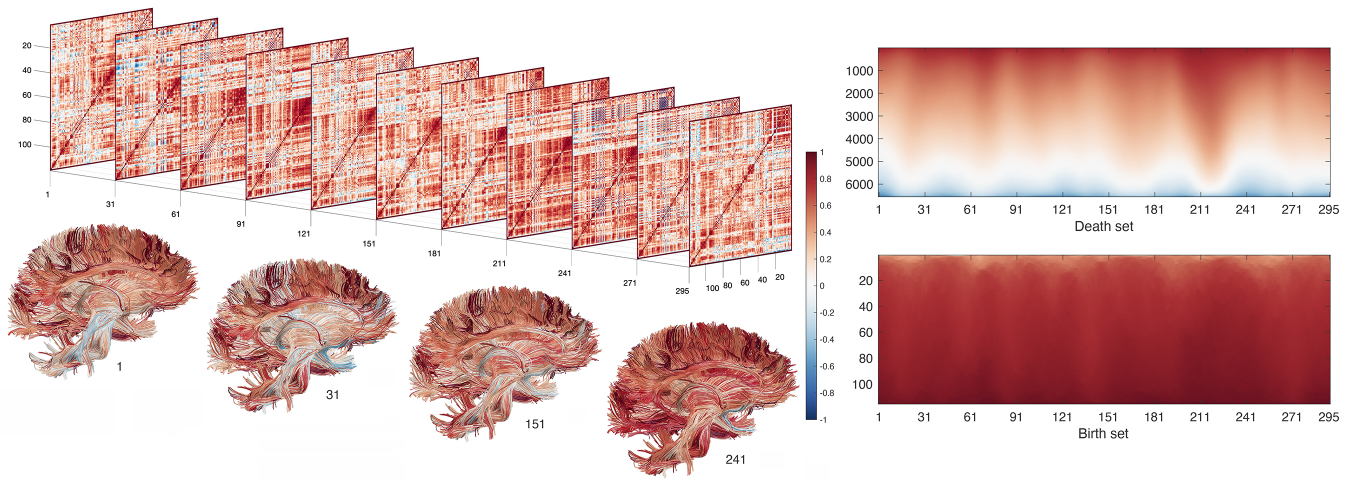}
\caption{The birth-death decomposition partitions the edge set into the birth and death edge sets. The birth set forms the maximum spanning tree (MST) and contains edges that create connected components (0D topology). The death set contains edges not belong to maximum spanning tree (MST) and destroys loops (1D topology).}
\label{fig:BDschematic}
\end{center}
\end{figure}

\begin{theorem}[Birth-death decomposition]
The edge weight set $W =  \{ w_{(1)}, \cdots, w_{(q)} \}$ has the unique decomposition
\bqn W = W_b \cup W_d,  \quad  W_b \cap W_d = \emptyset \label{eq:decompose} \eqn
where birth set $W_b = \{ b_{(1)}, b_{(2)}, \cdots, b_{(q_0)} \}$ is the collection of 0D sorted birth values and death set $W_d = \{ d_{(1)}, d_{(2)}, \cdots, d_{(q_1)} \}$ is the collection of 1D sorted death values with $q_0 = p-1$ and $q_1 = (p-1)(p-2)/2$. Further $W_b$ forms the 0D persistent diagram while $W_d$ forms the 1D persistent diagram. \label{thm:decompose}
\end{theorem}

In a complete graph with $p$ nodes, there are $q=p(p-1)/2$ unique edge weights. There are $q_0 = p-1$ number of edges that produces 0-cycles. This is equivalent to the number of edges in the maximum spanning tree of the graph. Thus, $q_1 = q - q_0 = \frac{(p-1)(p-2)}{2}$ number of edges destroy loops. The 0D persistent diagram  is given by $\{ (b_{(1)}, \infty),$ $\cdots,$ $(b_{(q_0)}, \infty) \}$. Ignoring $\infty$, $W_b$ is the 0D persistent digram. The 1D persistent diagram of the graph filtration is given by $\{ (-\infty, d_{(1)}),$ $\cdots, (-\infty, d_{(q_1)}) \}$. Ignoring $-\infty$, $W_d$ is the 1D persistent digram. We can show that the birth set is the maximum spanning tree (MST) (Figure \ref{fig:BDschematic}) \cite{song.2023}. \\

{\em Numerical implementation.} The identification of $W_b$ is based on the modification to Kruskal's or Prim's algorithm and identify the MST \cite{lee.2012.TMI}. Then $W_d$ is identified as $W / W_d$. Figure  \ref{fig:BDschematic} displays how the birth and death sets change over time for a single subject used in the study. Given edge weight matrix $W$ as input, the Matlab function {\tt WS\_decompose.m} outputs the birth set $W_b$ and the death set $W_d$.

\section{Topological Embedding}

Topological embedding provides a low-dimensional representation of complex network structures while preserving key topological features. Given a collection of brain networks represented as a set of connectivity matrices, we construct an embedding that captures the persistence of topological features across networks, allowing for a robust comparison of network dynamics.

\begin{figure}[t]
\begin{center}
\includegraphics[width=1\linewidth]{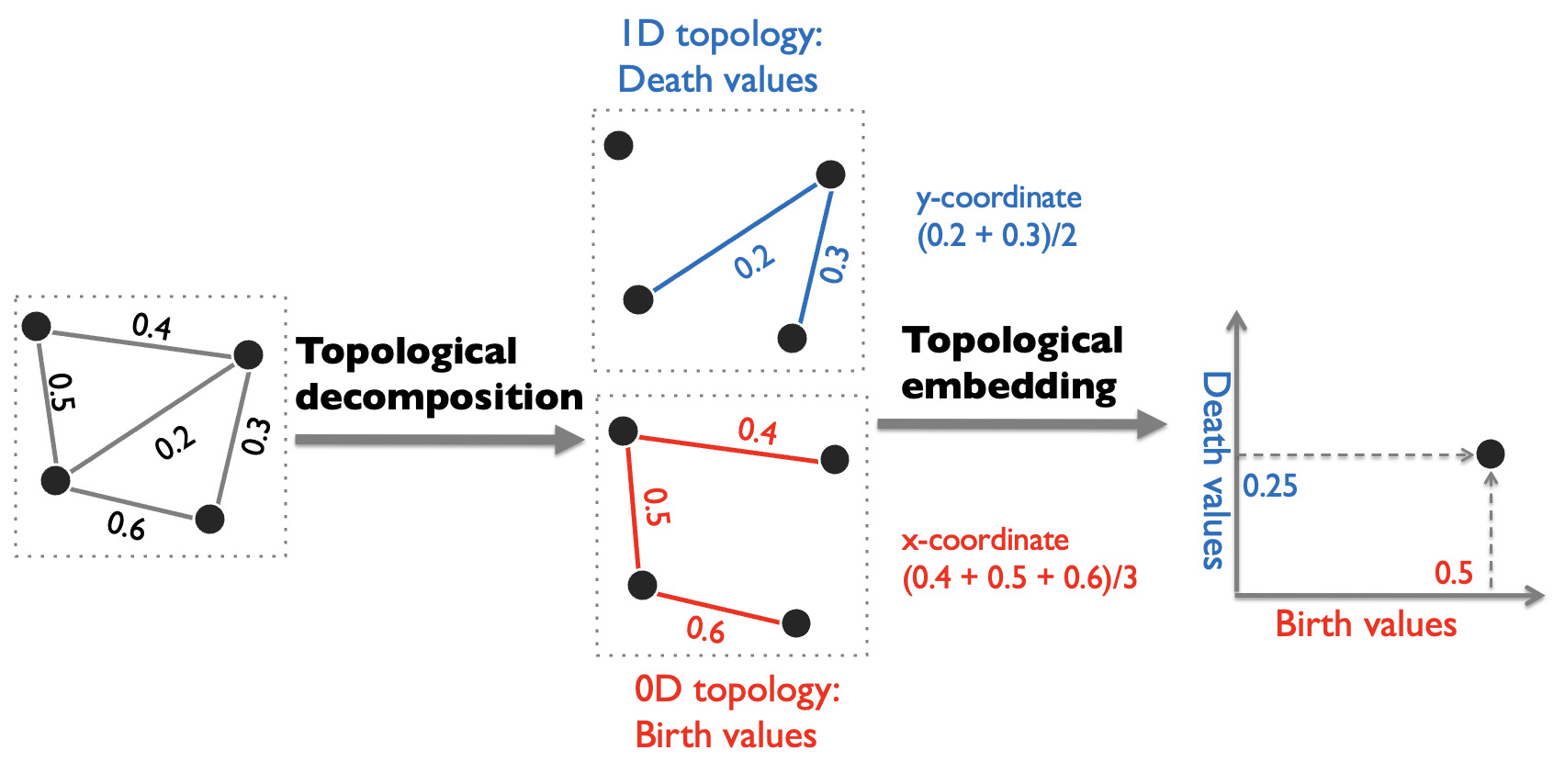}
\caption{The overall schematic of topological decomposition  followed by topological embedding. The birth values correspond to the maximum spanning tree (MST) of the underlying graph while the death values correspond to none MST edges. The process simplifies the weighted weighted network into a topological summary represented as a 2D point, where $x$-axis in TPD represents the average birth values, while the $y$-axis represents the average death values of the network.}
\label{fig:schematic}
\end{center}
\end{figure}

For a given collection of graphs or networks $\{\mathcal{X}_i\}_{i=1}^{N}$, we extract their topological features using persistent homology. The 0D homology captures the birth times of connected components, while the 1D homology identifies loops that form and disappear as edges are added. These birth and death values form a multiscale representation of network topology across the collection. To construct a meaningful embedding of these networks, we define the {\bf topological avearge}, which provides an aggregate representation of the collection by computing the mean birth and death values {\it within} networks \cite{chung.2023.NI}:

\begin{equation}
    b_i = \frac{1}{q_0} \sum_{j=1}^{q_0} b_{(j)}^i, \quad 
    d_i = \frac{1}{q_1} \sum_{j=1}^{q_1} d_{(j)}^i,
\end{equation}
where $b_{(j)}^i$ and $d_{(j)}^i$ are the $j$-th birth and death values extracted from network $\mathcal{X}_i$, and $q_0$ and $q_1$ represent the number of birth and death values in the graph filtration, respectively (Figure \ref{fig:schematic}). These values form a set of topological descriptors $(b_i, d_i)$ that summarize each network in a 2D topological space.

The set of topological descriptors $\{(b_i, d_i)\}_{i=1}^{N}$ forms an embedding of the network collection into a 2D topological space. The topological mean and standard deviation of the birth values across networks are defined as

\begin{equation}
    \mu_b = \frac{1}{N} \sum_{i=1}^{N} b_i, \quad
    \sigma_b = \sqrt{\frac{1}{N} \sum_{i=1}^{N} (b_i - \mu_b)^2}.
\end{equation}
Similarly, the mean and standard deviation of the death values are given by
\begin{equation}
    \mu_d = \frac{1}{N} \sum_{i=1}^{N} d_i, \quad
    \sigma_d = \sqrt{\frac{1}{N} \sum_{i=1}^{N} (d_i - \mu_d)^2}.
\end{equation}
Subsequently, we normalize the embedding points $(b_i, d_i)$ as follows:
\begin{equation}
    \tilde{b}_i = \frac{b_i - \mu_b}{\sigma_b}, \quad 
    \tilde{d}_i = \frac{d_i - \mu_d}{\sigma_d},
\end{equation}
where $\mu_b$ and $\sigma_b$ are the topological group mean and standard deviation of birth values across networks, and $\mu_d$ and $\sigma_d$ are the mean and standard deviation of death values. This normalization ensures that variations due to differences in scale or network density are accounted for, making the topological embedding more robust to inter-subject variability.

\begin{figure}
\centering
\includegraphics[width=0.7\linewidth]{TPD-3statesv2.jpg} 
\caption{The topological embedding  showing the evolution of rs-fMRI network changes over time for a single subject. The horizontal axis represents changes of 0D topology (connected components) while the vertical axis represents changes of 1D topology (loops).}
\label{fig:pipeline}
\vspace{-.25cm}
\end{figure}

The topological embedding provides an aggregate representation of the network collection, allowing for comparisons across groups, conditions, or cognitive states. Networks that are topologically similar will cluster together, while deviations from the centroid indicate structural or functional differences. The distance of an individual network’s embedding $(\tilde{b}_i, \tilde{d}_i)$ from the group eman $(\mu_b, \mu_d)$ provides a measure of {\it topological deviation}, which can be used to assess network stability 
\begin{equation}
    \delta_i = \sqrt{(\tilde{b}_i - \mu_b)^2 + (\tilde{d}_i - \mu_d)^2}.
\end{equation}
A lower $\delta_i$ indicates that a network closely follows the overall topological structure of the population, while a higher $\delta_i$ suggests significant deviation, which may correspond to atypical functional brain connectivity patterns.

The topological embedding framework enables a compact yet interpretable representation of time-evolving or subject-specific connectivity networks as well. It can be used to compare brain networks across populations by clustering networks based on their topological proximity. Analyze cognitive variability by examining how networks deviate from the topological centroid in different cognitive states. Track dynamic changes in functional connectivity by embedding time-varying networks and identifying stable and transient topological structures. Figure \ref{fig:pipeline} illustrates an example of a topological embedding, where variations in birth and death values reveal dynamic changes in brain network configurations over time.

\section{Topological Distance}

\begin{figure}[t]
\begin{center}
\includegraphics[width=1\linewidth]{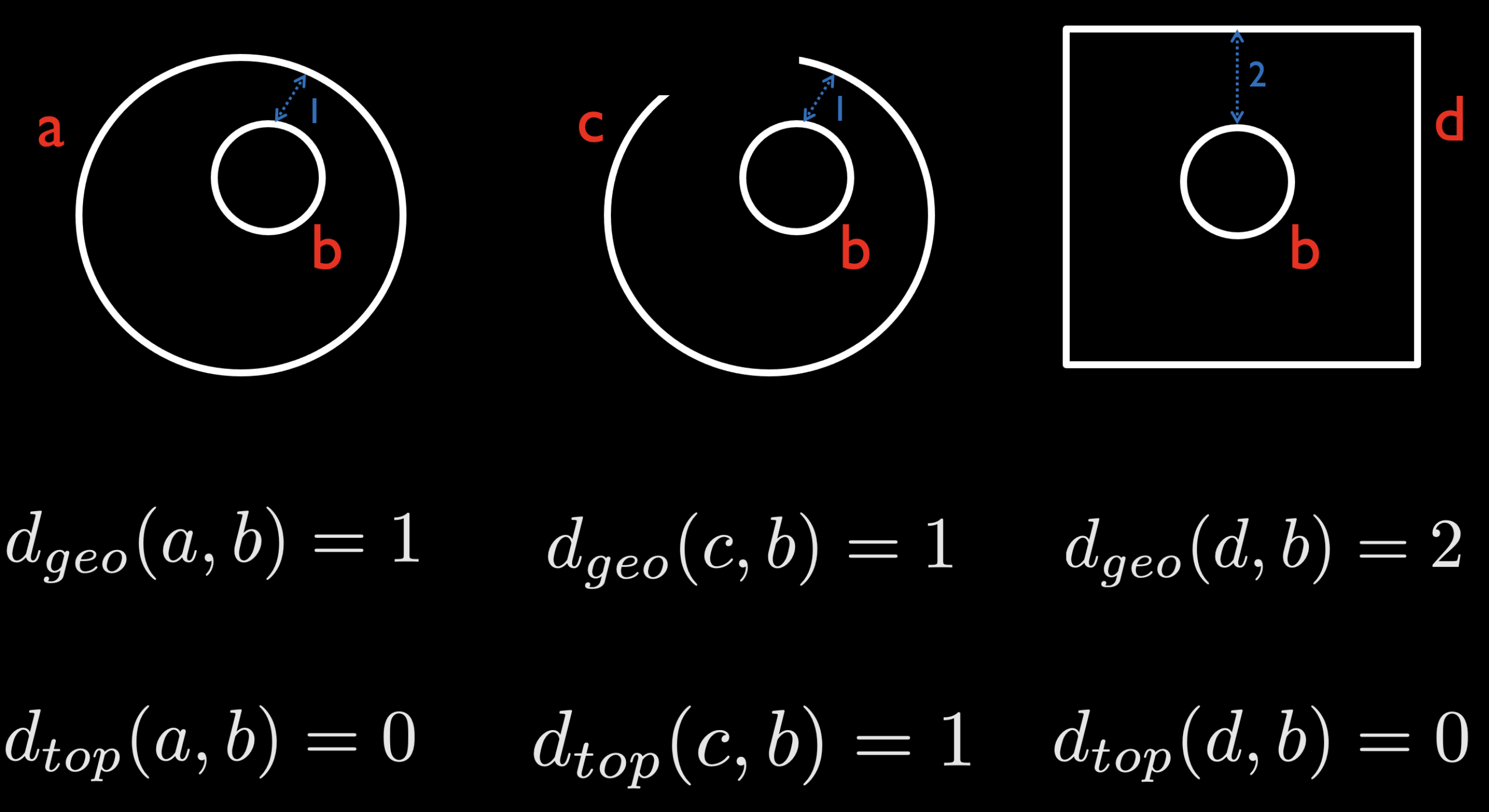}
\caption{Comparison between geometric distance $d_{geo}$ and topological distance $d_{top}$. We used the shortest distance between objects as the geometric distance. The left and middle objects are topologically different while the left and right objects are topologically equivalent. The geometric distance cannot discriminate topologically different objects (left and middle) and produces false negatives. The geometric distance incorrectly discriminate topologically equivalent objects  (left and right)and produces false positives.}
\label{fig:geotop}
\end{center}
\end{figure}

The main difference between the geometric and topological distance is if the distance can discriminate the presence of topological difference and not able to discriminate the presence of topological indifference (Figure \ref{fig:geotop}). This can be achieved using the Wasserstein distance between persistent diagrams.

\subsection{Wasserstein distance}
Given two probability distributions $X \sim f_1$ and $Y \sim f_2$, the $r$-{\em Wasserstein distance} $D_W$, which is the probabilistic version of the optimal transport,  is defined as
$$\mathcal{D} (f_1, f_2) =  \Big( \inf \mathbb{E} |  X - Y |^r \Big)^{1/r},$$
where the infimum is taken over every possible joint distributions of $X$ and $Y$.  The Wasserstein distance
is the optimal expected cost of transporting points generated from $f_1$ to those generated from $f_2$ \cite{canas.2012}. There are numerous distances and similarity measures defined between probability distributions such as the Kullback-Leibler (KL) divergence and the mutual information \cite{kullback.1951}. While the Wasserstein distance is a metric satisfying positive definiteness, symmetry,
and triangle inequality, the KL-divergence and the mutual information are not metric. Although they are easy to compute, the biggest limitation of the KL-divergence and the mutual information is that the two probability distributions has to be defined on the same sample space. If the two distributions do not have the same support, it may be difficult to even define them. If $f_1$ is discrete while $f_2$ is continuous, it is difficult to define them. On the other hand, the Wasserstein distance can be computed for any arbitrary distributions that may not have the common sample space making it extremely versatile. 

The Wasserstein distance is often used distance for measuring the discrepancy between persistent diagrams. Consider persistent diagrams $P_1$ and $P_2$ given by  
$$P_1: x_1 = (b_1^1, d_1^1), \cdots, x_q = (b_q^1,d_q^1),  \quad P_2: y_1 = (b_1^2, d_1^2), \cdots, y_q = (b_q^2, d_q^2).$$ 
Their empirical distributions are given in terms of Dirac-delta functions 
$$f_1 (x) = \frac{1}{q} \sum_{i=1}^q \delta (x-x_i), \quad f_2(y) = \frac{1}{q} \sum_{i=1}^q \delta (y-y_i).$$
Then we can show that the $r$-{\em Wasserstein distance} on persistent diagrams  is given by
\bqn \mathcal{D} (P_1, P_2) = \inf_{\psi: P_1 \to P_2} \Big( \sum_{x \in P_1} \| x - \psi(x) \|^r \Big)^{1/r} \label{eq:D_Winf} \eqn
over every possible bijection $\psi$, which is permutation, between $P_1$ and $P_2$ \cite{vallender.1974,canas.2012,berwald.2018}. 
Optimization (\ref{eq:D_Winf}) is the standard assignment problem, which is  usually solved by Hungarian algorithm in $\mathcal{O} (q^3)$ \cite{edmonds.1972}. However, for graph filtration, the distance can be computed {\em exactly} in $\mathcal{O}(q \log q)$  by simply matching the order statistics on the birth or death values \cite{rabin.2011,song.2023,song.2021.MICCAI}. Note the 0D persistent diagram for graph filtration is just 1D scatter points of birth values while 
the 1D persistent diagram for graph filtration is just 1D scatter points of death values. Thus, the Wasserstein distance can be simplified as follows.  

\begin{theorem} \cite{song.2023} The r-Wasserstein distance between the 0D persistent diagrams for graph filtration is given by 
$$\mathcal{D}_0 (P_1, P_2) = \Big[ \sum_{i=1}^{q_0} (b_{(i)}^1 - b_{(i)}^2)^r \Big]^{1/r},$$
where $b_{(i)}^j$ is the $i$-th smallest birth values in persistent diagram $P_j$.  The 2-Wasserstein distance between the 1D persistent diagrams for graph filtration is given by 
$$\mathcal{D}_1 (P_1, P_2) =  \Big[ \sum_{i=1}^{q_1} (d_{(i)}^1 - d_{(i)}^2)^r  \Big]^{1/r},$$
where $d_{(i)}^j$ is the $i$-th smallest death values in persistent diagram $P_j$. 
\label{theorem:optimal}
\end{theorem}

Here, we will simply use the combined 0D and 1D topological distance $\mathcal{D} = \mathcal{D}_0 + \mathcal{D}_1$ for inference and clustering. For collection of graphs {\tt con\_i} and {\tt con\_j}, the pairwise Wasserstein distance between graphs is computed as 
\begin{verbatim} 
lossMtx = WS_pdist2(con_i,con_j)

  struct with fields:

     D0: [10×10 double]
     D1: [10×10 double]
    D01: [10×10 double]
\end{verbatim}
Each entry of {\tt lossMtx} stores 0D, 1D and combined topological distances.

\subsection{Topological registration}
The Wasserstein distance allows for one-to-one correspondence of edges between two networks. Thus, it can be utilized for network registration in $\mathcal{O}(q \log q)$ runtime. One limitation occurs in matching networks of different sizes.  When networks have different sizes in terms of the number of nodes, we can still find the optimal matching through data augmentation on smaller networks. Similar data augmentations are performed in matching various topological features of different sizes \cite{hu.2019}, matching trees of different sizes \cite{guo.2020}, or matching point sets of different sizes \cite{song.2023}.

\begin{figure}[t]
\begin{center}
\includegraphics[width=1\linewidth]{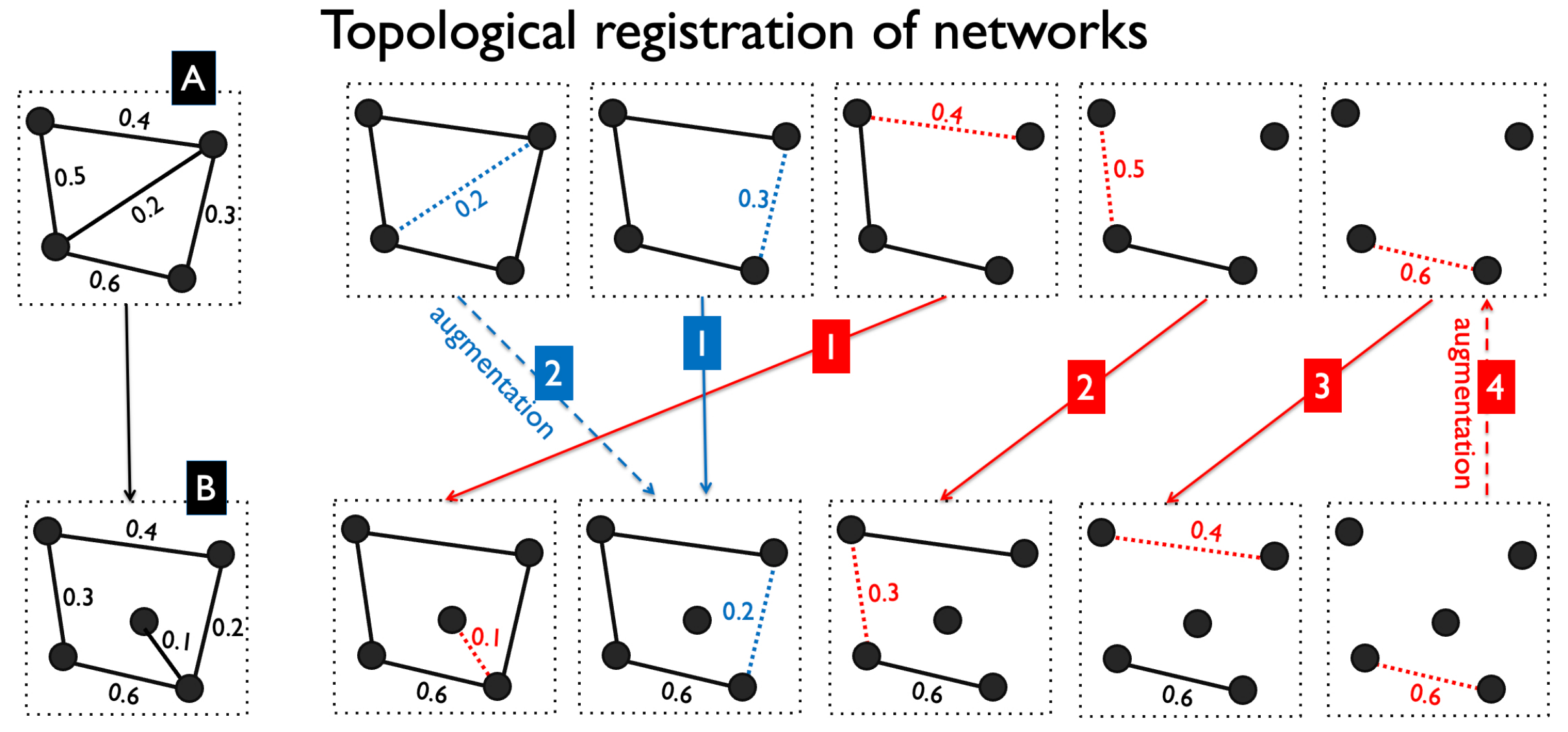}
\caption{Schematic of registering network A with 4 nodes to network B with 5 nodes. Through the birth-death decomposition, we can match edges belonging to 0D topology (red) and 1D topology (blue) separately. Numbers inside blue and red squares indicate the order of matching. Any unmatched long-lived features are matched to the largest persistent feature, and short-lived features are matched to the smallest persistent feature through augmentations.}
\label{fig:aug}
\end{center}
\end{figure}

Connected components that are born later (with larger birth values) have shorter persistences (lifetimes), which are considered as smaller signals \cite{ghrist.2008,clough.2019,hu.2019}. Thus, we need to map any unmatched, short-lived connected components  to the least persistent connected component. Thus, any remaining unmatched large birth values in one network are mapped to the largest edge weight in the other graph (step 4 in red square in Figure \ref{fig:aug}). Unmatched death values are handled analogously: the first largest unmatched death values are mapped to the largest death value in the other graph, while any unmatched small death values are mapped to the smallest edge weight in the other graph (step 2 in blue square in Figure \ref{fig:aug}). By augmenting the networks in this manner, we ensure a consistent and meaningful comparison of topological features across networks of different sizes.

\subsection{Topological averaging}
\label{sec:TopoAveraging}

The concept of topological registration of networks can be used to topologically average networks. For a collection of networks represented as a set of graphs \(\{\mathcal{X}_i\}_{i=1}^n\), we introduce an approach to identify the topological mean of these networks. This concept is akin to the Wasserstein barycenter \cite{agueh.2011,cuturi.2014} and the Fréchet mean \cite{le.2000,turner.2014,zemel.2019,dubey.2019} but offers scalable \(\mathcal{O}(q \log q)\) computation for even large networks, where $q$ is the number of edges \cite{chung.2023.NI}. Unlike conventional averaging methods that focus on element-wise arithmetic means of edge weights, our method leverages the Wasserstein distance to capture the underlying topological features of the networks.

We define the \textit{topological addition} $\oplus$ to amalgamate two graphs, \(\mathcal{X}_1 = (V, w^1)\) and \(\mathcal{X}_2 = (V, w^2)\)
, into \(\mathcal{X}_1 \oplus \mathcal{X}_2 = (V, w)\) via topological registration. The birth-death decomposition of \(\mathcal{X}_1 \oplus \mathcal{X}_2\) is given by element-wise addition of the sorted birth and death sets:
\bq
B_1 + B_2 &=& \{ b_{(1)}^1 + b_{(1)}^2, \cdots, b_{(q_0)}^1 + b_{(q_0)}^2 \}\\
D_1 + D_2 &=& \{ d_{(1)}^1 + d_{(1)}^2, \cdots, d_{(q_1)}^1 + d_{(q_1)}^2 \}.
\eq

\begin{theorem}
\label{theorem:sum}
The topological mean
\[
\mu = \frac{1}{n}\oplus_{j=1}^n \mathcal{X}_j
\]
of graphs \(\mathcal{X}_1, \mathcal{X}_2, \ldots, \mathcal{X}_n\)
is the minimizer of the sum of topological distances to all other graphs
$$
\mu = \arg \min_{Y} \sum_{i=1}^n d^2(Y, \mathcal{X}_i).
\label{eq:centroid_relaxed}
$$
\end{theorem}

{\it Proof.} 
We express the topological distance \(d^2(Y, \mathcal{X}_i)\) for a candidate graph \(Y\):
\[
d^2(Y, \mathcal{X}_i) = \sum_{j=1}^{q_0} (b_j^Y - b_j^i)^2 + \sum_{j=1}^{q_1} (d_j^Y - d_j^i)^2.
\]
Since the cost function is a linear combination of quadratic functions, the global minimum exists and is unique. The sum is minimized when
\[
\frac{\partial}{\partial b_j^Y} \sum_{i=1}^n d^2(Y, \mathcal{X}_i) = 
\frac{\partial}{\partial d_j^Y} \sum_{i=1}^n d^2(Y, \mathcal{X}_i) = 0.
\]
Solving for \(b_j^Y\) and \(d_j^Y\), we get
\[
b_j^Y = \frac{1}{n} \sum_{i=1}^n b_j^i, \quad d_j^Y = \frac{1}{n} \sum_{i=1}^n d_j^i.
\]
Thus, the birth and death sets of the candidate graph $Y$ are
\[
B = \{ \frac{1}{n} \sum_{i=1}^n b_j^i \}_{j=1}^{q_0}, \quad D =  \{ \frac{1}{n} \sum_{i=1}^n d_j^i \}_{j=1}^{q_1}.
\]
Therefore, the topological mean \(\mu\) defined as \( (\mathcal{X}_1 \oplus \mathcal{X}_2 \oplus \ldots \oplus \mathcal{X}_n)/n\) is the candidate graph. \hfill $\square$

Similar to the Fréchet mean \cite{le.2000,turner.2014,zemel.2019,dubey.2019}, which may not be unique in certain configurations, we emphasize that the topological mean may also exhibit non-uniqueness under specific conditions \cite{chung.2023.NI}

\section{Topological Inference}
There are a few studies that used the Wasserstein distance \cite{mi.2018,yang.2020}. The existing methods are mainly applied to geometric data without topological consideration. It is not obvious how to apply the method to  perform statistical inference for a population study. We will present a new statistical inference procedure for testing the topological inference of two groups, the usual setting in brain network studies. 
Consider a collection of graphs $\mathcal{X}_1, \cdots, \mathcal{X}_n$  that are grouped into two groups 
$C_1$ and $C_2$ such that 
$$ C_1 \cup C_2 = \{\mathcal{X}_1, \cdots, \mathcal{X}_n\} , \quad C_1 \cap C_2 = \emptyset.$$
We assume there are $n_i$ graphs in $C_i$ and $n_1  + n_2 = n$. In the usual statistical inference, we are interested in testing the null hypothesis of the equivalence of topological summary $\mathcal{T}$:
$$H_0: \mathcal{T}(C_1)  = \mathcal{T}(C_2).$$
Under the null, there are ${n \choose n_1}$ number of permutations to permute $n$ graphs into two groups, which is an extremely large number and most computing systems including MATLAB/R cannot compute them exactly if the sample size is larger than 50 in each group. If $n_1 = n_2$, the total number of permutations is given asymptotically by Stirling's formula \cite{feller.2008}
$${n \choose n_1} \sim \frac{4^{n_1}}{\sqrt{\pi n_1}}.$$
The number of permutations {\em exponentially} increases as the sample size increases, and thus it is impractical to generate every possible permutation. In practice, up to hundreds of thousands of random permutations are generated using the uniform distribution on  the permutation group with probability $1/{n \choose n_1}$. The computational bottleneck in the permutation test is mainly caused by the need to recompute the test statistic for each permutation. This usually cause a serious computational bottleneck when we have to recompute the test statistic for large samples for  more than million permutations. We propose a more scalable approach. 

 \begin{figure}[t!]
\centering
\includegraphics[width=1\linewidth]{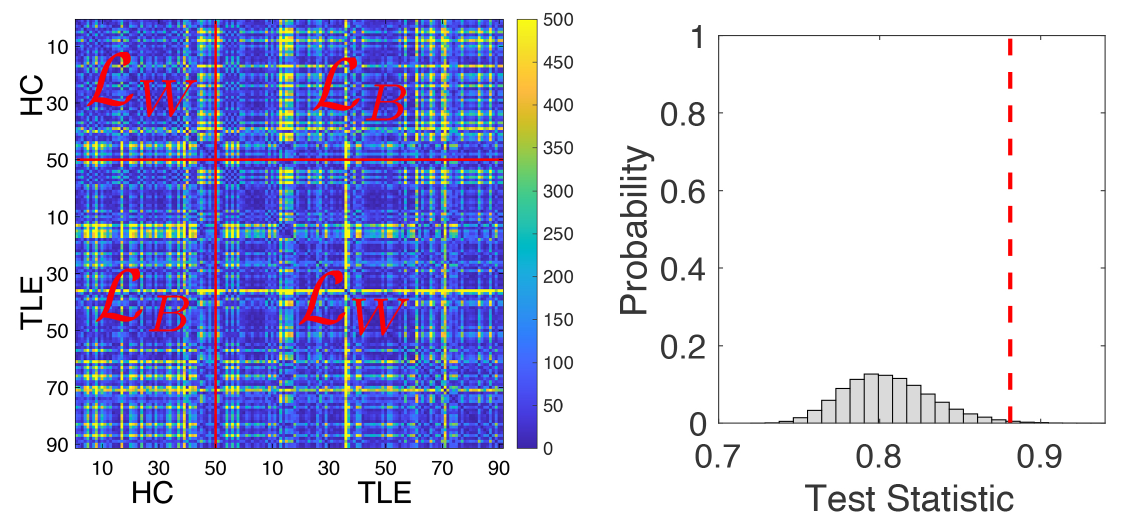}
\caption{Pairwise Wasserstein distance between 50 healthy controls (HC) and 101 temporal lobe epilepsy (TLE) patients. There are subtle pattern difference in the off-diagonal patterns (between group distances $\mathcal{L}_B$) compared to diagonal patterns (within group distances $\mathcal{L}_W$).  The permutation test with 100 million permutations was used to determine the statistical significance using the ratio statistic. The  red line is the observed ratio. The histogram is the empirical null distribution obtained from the permutation test.}
\label{fig:exampletheroem}
\end{figure}

Define  the within-group distance  $\mathcal{L}_{W}$ as
\[
2 \mathcal{L}_{W} =   \sum_{\mathcal{X}_i, \mathcal{X}_j \in C_1} \mathcal{D}(\mathcal{X}_i,\mathcal{X}_j)  + \sum_{\mathcal{X}_i, \mathcal{X}_ j \in C_2} \mathcal{D}(\mathcal{X}_i,\mathcal{X}_j). 
 \]
 The within-group distance corresponds to the sum of all the pairwise distances in the block diagonal matrices in Figure \ref{fig:exampletheroem}. The average within-group distance is then given by
$$ \overline{\mathcal{L}}_{W} = \frac{\mathcal{L}_{W}}{n_1(n_1-1) + n_2 (n_2-1)}.$$ 
 The between-group distance $\mathcal{L}_{B}$ is defined as
 \[ 2 \mathcal{L}_{B}=  \sum_{\mathcal{X}_i \in C_1} \sum_{\mathcal{X}_j \in C_2} \mathcal{D}(\mathcal{X}_i,\mathcal{X}_j) + \sum_{\mathcal{X}_i \in C_2} \sum_{\mathcal{X}_j \in C_1} \mathcal{D}(\mathcal{X}_i,\mathcal{X}_j) .\]
 The between-group distance corresponds to the off-diaognal block matrices in Figure \ref{fig:exampletheroem}. 
 The average between-group distance is then given by
$$ \overline{\mathcal{L}}_{B} = \frac{\mathcal{L}_{B}}{n_1 n_2}.$$

 Note that the sum of  within-group and  between-group distance is the sum of all the pairwise distances in Figure \ref{fig:exampletheroem}:
 $$ 2\mathcal{L}_{W} + 2 \mathcal{L}_{B} = \sum_{i=1}^n \sum_{j=1}^n  \mathcal{D}(\mathcal{X}_i,\mathcal{X}_j).$$ 
When we permute the group labels, the total sum of all the pairwise distances do not change and fixed. If the group difference is large, the between-group distance $\mathcal{L}_{B}$ will be large and the within-group distance $\mathcal{L}_{W}$ will be small. Thus, to  measure the disparity between groups as the  ratio  \cite{song.2023}
\[ \phi_\mathcal{L} =  \frac{\mathcal{L}_{B}}{\mathcal{L}_{W}} .\]
The ratio statistic is related to the elbow method in clustering and behaves like traditional $F$-statistic, which is the ratio of squared variability of model fits. If $\phi_\mathcal{L}$ is large, the groups differ significantly in network topology. If $\phi_\mathcal{L}$ is small, it is likely that there is no group differences.

Since the ratio is always positive, its probability distribution cannot be Gaussian. Since the distributions of the ratio $\phi_\mathcal{L}$ is unknown, the permutation test can be used to determine the empirical distributions. Figure \ref{fig:exampletheroem}-right displays the empirical  distribution of $\phi_\mathcal{L}$. The $p$-value is the area of the right tail thresholded by the observed ratio $\phi_\mathcal{L}$ (dotted red line) in the empirical distribution. Since we only compute the pairwise distances only once and only shuffle each entry over permutations. This is equivalent  to rearranging rows and columns of entries corresponding to the permutations in Figure \ref{fig:exampletheroem}. The simple rearranging of rows and columns of entries and sum them in the block-wise fashion should be  faster than the usual two-sample $t$ test which has to be recomputed for each permutation.

To speed up the permutation further, we adapted the transposition test, the online version of permutation test \cite{chung.2019.CNI}. In the transposition test, we only need to work out how 
$\mathcal{L}_{B}$ and $\mathcal{L}_{W}$ changes over a transposition, a permutation that only swaps one entry from each group. When we  transpose $k$-th and $l$-th graphs between the groups (denoted as $\tau_{kl}$), all the $k$-th and $i$-th rows and columns will be swapped. The within-group distance after the transposition $\tau_{kl}$ is given by
$$\tau_{kl} (\mathcal{L}_{W})  = \mathcal{L}_{W} + \Delta_{W},$$
where $\Delta_{W}$ is the terms in the $k$-th and $i$-th rows and columns that are required to swapped. We only need to swap up to $\mathcal{O}(2n)$ entries while the standard permutation test that requires the computation over $\mathcal{O}(n^2)$ entries. Similarly we have incremental changes 
$$\tau_{kl} (\mathcal{L}_{B})  = \mathcal{L}_{B} + \Delta_{B}.$$
The ratio statistic over the transposition is then sequentially updated over random transpositions. To further accelerate the convergence and avoid potential bias, we introduce one permutation to the sequence of 1000 consecutive transpositions.

\begin{figure}[t]
\centering
\includegraphics[width=1\linewidth]{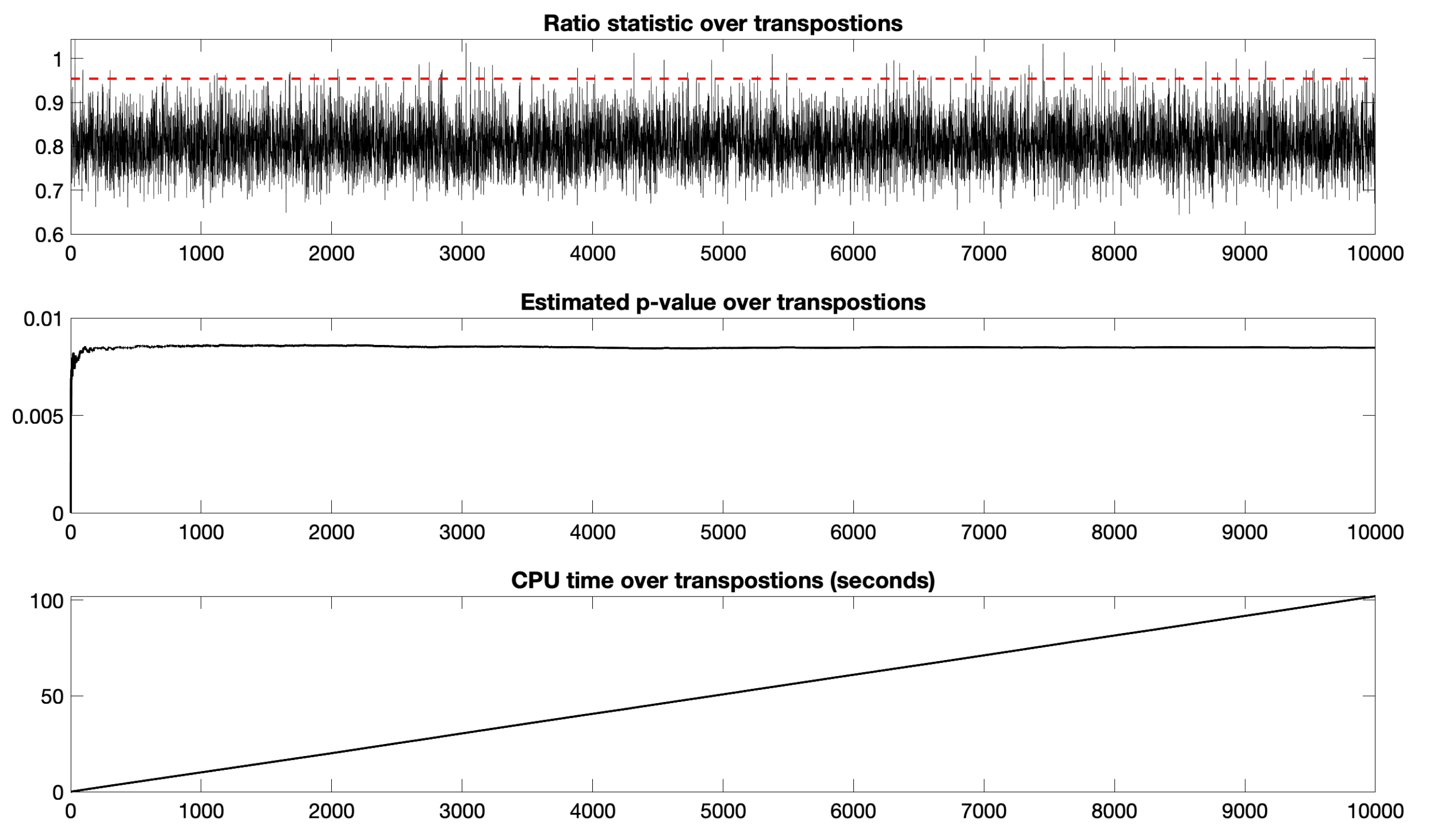}
\caption{The plot of ratio statistic $\phi_{\mathcal{L}}$ (top) over 100 million transpositions in testing the topological difference between HC and TLE. The plot is only shown at every 10000 transposition. The redline is the observed ratio static 0.9541. The estimated $p$-value (middle) converges to 0.0086 after 100 million transpositions. The CPU time (bottom) is linear and takes 102 seconds for 100 million transpositions.}
\label{fig:convergence}
\end{figure}

The observed ratio statistic is computed using {\tt WS\_ratio.m}, which inputs the distance matrix {\tt lossMtx}, sample size in each group. The whole procedure for performing the transposition test is  implemented as {\tt WS\_transpositions.m} and takes less than  one second in a desktop computer for million permutations. The function inputs the distance matrix {\tt lossMtx}, sample size in each group, number of transpositions and the number of permutations that are interjected into transpositions. Figure \ref{fig:convergence} displays the convergence plot of the transposition test.

\section{Topological Clustering}

We validate  the proposed topological distances in simulations with the ground truth in a clustering setting. The Wasserstein distance was previously used for clustering  for {\em geometric objects} without topology in $\mathcal{O}(q^3)$ \cite{mi.2018,yang.2020}. The proposed topological method builds the Wasserstein distances on persistent diagrams in $\mathcal{O}(q \log q)$ making our method scalable. Consider a collection of graphs $\mathcal{X}_1, \cdots, \mathcal{X}_n$ that will be clustered into $k$ clusters $C=(C_1, \cdots, C_k)$. Let $\mu_j = \mathbb{E} C_j$ be the topological mean of $C_j$ computing using the Wasserstein distance.  Let $\mu = (\mu_1, \cdots, \mu_k)$ be the cluster mean vector.
The within-cluster Wasserstein distance  is given by
\bq l_W (C; \mu) =  \sum_{j=1}^k \sum_{X \in C_j} \mathcal{D}(X, \mu_j) =  \sum_{j=1}^k |C_j| \mathbb{V} C_j
\label{eq:l_W}
 \eq
with the topological variance $\mathbb{V} C_j$ of cluster $C_j$. 
The within-cluster Wasserstein distance generalizes the within-group distance defined on two groups to $k$ number of groups (or clusters). When $k=2$, we have $l_W (C; \mu)  = 2\mathcal{L}_W$.

The topological clustering through the Wasserstein distance is then performed by minimizing $l_W(C)$  over every possible  $C$. The Wasserstein graph clustering algorithm can be implemented as the two-step optimization often used in variational inferences \cite{bishop.2006}. The algorithm follows the proof below.

\begin{theorem} Topological clustering with the Wasserstein distance converges locally.
\end{theorem}
\begin{proof}  1) Expectation step: Assume $C$ is estimated from the previous iteration. In the current iteration, the cluster mean $\mu$ corresponding to $C$ is updated as 
$\mu_j \leftarrow \mathbb{E} C_j$
for each $j$. The cluster mean gives the lowest bound on distance $l_W (C; \nu)$ for any $\nu = (\nu_1, \cdots, \nu_k)$:
\bqn l_W (C; \mu) =  \sum_{j=1}^k \sum_{X \in C_j} \mathcal{D}(X, \mu_j) \leq  \sum_{j=1}^k \sum_{X \in C_j} \mathcal{D}(X, \nu_j) = l_W(C; \nu). \label{eq:EMexp} \eqn
2) We check if the cluster mean $\mu$ is changed from the previous iteration. If not, the algorithm simply stops. Thus we can force $l_W (C; \nu)$ to be strictly decreasing over each iteration. 3) Minimization step: The clusters are updated from $C$ to $C' = (C'_{J_1}, \cdots, C'_{J_k})$
by reassigning each graph $\mathcal{X}_i$ to the closest cluster $C_{J_i}$ satisfying
$J_i =\arg \min_j \mathcal{D}(\mathcal{X}_i, \mu_j).$
Subsequently, we have
\bqn l_W (C'; \mu) =  \sum_{J_i=1}^k \sum_{X \in C_{J_i}'} \mathcal{D}( X, \mu_{J_i}) \leq  \sum_{j=1}^k \sum_{X \in C_j} \mathcal{D}(X, \mu_j) = l_W (C; \mu). \label{eq:EMmin}\eqn
From (\ref{eq:EMexp}) and  (\ref{eq:EMmin}),  $l_W (C; \mu)$ strictly decreases over iterations. Any bounded strictly decreasing sequence converges.  
\end{proof}
Just like $k$-means clustering  that converges only to local minimum, there is no guarantee the Wasserstein graph clustering converges to the global minimum \cite{huang.2020.NM}. This is remedied by repeating the algorithm multiple times with different random seeds and identifying the cluster that gives the minimum over all possible seeds.

Let $y_i$ be the true cluster label for the $i$-th data. Let $\widehat y_i$ be the estimate of $y_i$ we determined from Wasserstein graph clustering. Let $y=(y_1, \cdots, y_n)$ and $\widehat y=(\widehat y_1, \cdots, \widehat y_n)$. In clustering, there is no direct association between true clustering labels and predicted cluster labels. Given $k$ clusters $C_1, \cdots, C_k$, its permutation $\pi(C_1), \cdots,$ $\pi(C_k)$ is also a valid cluster for $\pi \in \mathbb{S}_k$,  the permutation group of order $k$. There are $k!$ possible permutations in $\mathbb{S}_k$ \cite{chung.2019.CNI}. The clustering accuracy $A(y, \widehat y)$ is then given by
$$A(\widehat y,y) = \frac{1}{n} \max_{\pi \in \mathbb{S}_k} \sum_{i=1}^n \mathbf{1} ( \pi( \widehat y) = y).$$
This a modification to an assignment problem and can be solved using the Hungarian algorithm in $\mathcal{O}(k^3)$ run time \cite{edmonds.1972}. In Matlab, it can be solved using {\tt confusionmat.m}, which tabulates misclustering errors between the true cluster labels and predicted cluster labels. 
Let $F(\widehat y,y) = (f_{ij})$ be  the confusion matrix of size $k \times k$ tabulating the correct number of clustering in each cluster. The diagonal entries show the correct number of clustering while the off-diagonal entries show the incorrect number of clusters. 
To compute the clustering accuracy, we need to sum the diagonal entries. 
Under the permutation of cluster labels, we can get different confusion matrices. For large $k$, it is prohibitive expensive to search for all permutations. Thus we need to maximize the sum of diagonals of the confusion matrix under permutation: 
\bqn \frac{1}{n} \max_{Q \in \mathbb{S}_k} \mbox{tr} (QC) = \frac{1}{n} \max_{Q \in \mathbb{S}_k} \sum_{i,j} q_{ij} f_{ij}, \label{eq:optimization} \eqn
where $Q=(q_{ij})$ is the permutation matrix consisting of entries 0 and 1 such that there is exactly single 1 in each row and each column. This is a linear sum assignment problem (LSAP), a special case of linear assignment problem \cite{duff.2001,lee.2018.deep}. The clustering accuracy is computed using
\begin{verbatim}
function [accuracy C]=cluster_accuracy(ytrue,ypred)
\end{verbatim}
where {\tt ytrue} is the true cluster labels and {\tt ypred} is the predicted cluster labels. {\tt accuracy} is the clustering accuracy and {\tt C} is the confusion matrix \cite{duff.2001}. 

\begin{figure}[t]
\centering
\includegraphics[width=1\linewidth]{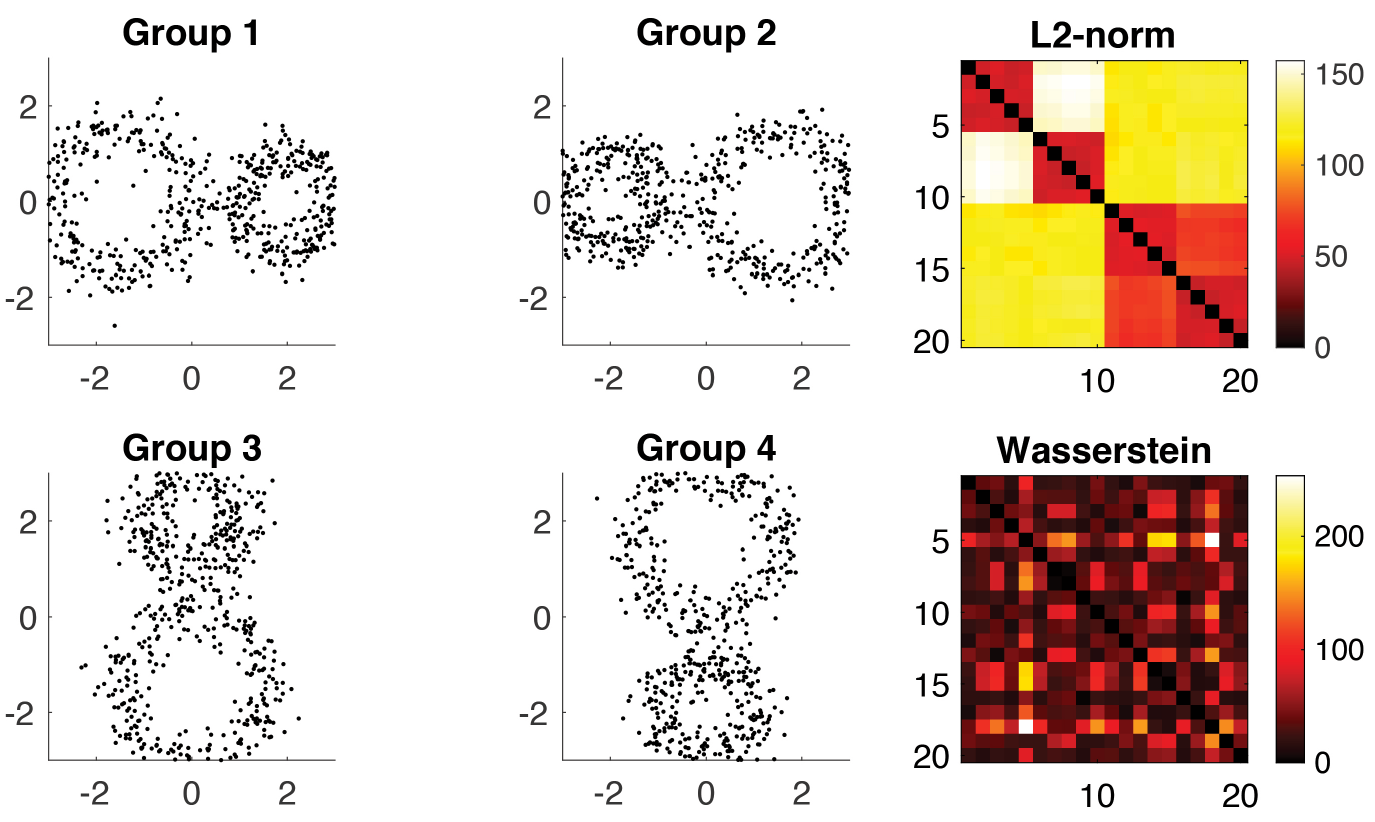}
\caption{Simulation study on  topological equivalence. The correct clustering method should {\em not} be able to cluster them since they are all topologically equivalent. Right: the pairwise Euclidean distance ($L_2$-norm) is used in $k$-means and hierarchical clustering. The Wasserstein distance is used in topological clustering.}
\label{fig:simulationnodiff}
\end{figure}

\begin{example}
We replace the Euclidean distance ($L_2$-norm) in $k$-means clustering with  the topological distance $\mathcal{D}$ and compared the performance with the traditional $k$-means clustering and hierarchical clustering \cite{lee.2011.MICCAI}. We generated 4 circular patterns of identical topology (Figure \ref{fig:simulationnodiff}) and different topology (Figure \ref{fig:simulationdiff}). Along the circles, we uniformly sampled 60 nodes and added  Gaussian noise $N(0, 0.3^2)$ on the coordinates. We generated 5 random networks per group. The Euclidean distance ($L_2$-norm)  between randomly generated points are used to build connectivity matrices for $k$-means  and hierarchical clustering. Figures \ref{fig:simulationnodiff} and \ref{fig:simulationdiff} shows the superposition of nodes from 20 networks.   For $k$-means and Wasserstein graph clustering, the average result of 100 random seeds are reported.

We tested for false positives when there is no topology difference in Figure \ref{fig:simulationnodiff}, where all the groups are simply obtained from Group 1 by rotations. All the groups are topologically equivalent and thus we should not detect any topological difference. Any detected signals are all false positives. 
The $k$-means had $0.90 \pm 0.15$ while the hierarchical clustering had perfect $1.00$ accuracy. Existing clusterng methods based on Euclidean distance are reporting significant false positives and should not be used in topological clustering task had the accuracy. On the other hand, the Wasserstein graph clustering had low $0.53 \pm 0.08$ accuracy.  We conclude that Wasserstein graph clustering are not reporting topological false positive like $k$-means and hierarchical clusterings. 

\begin{figure}[t]
\centering
\includegraphics[width=1\linewidth]{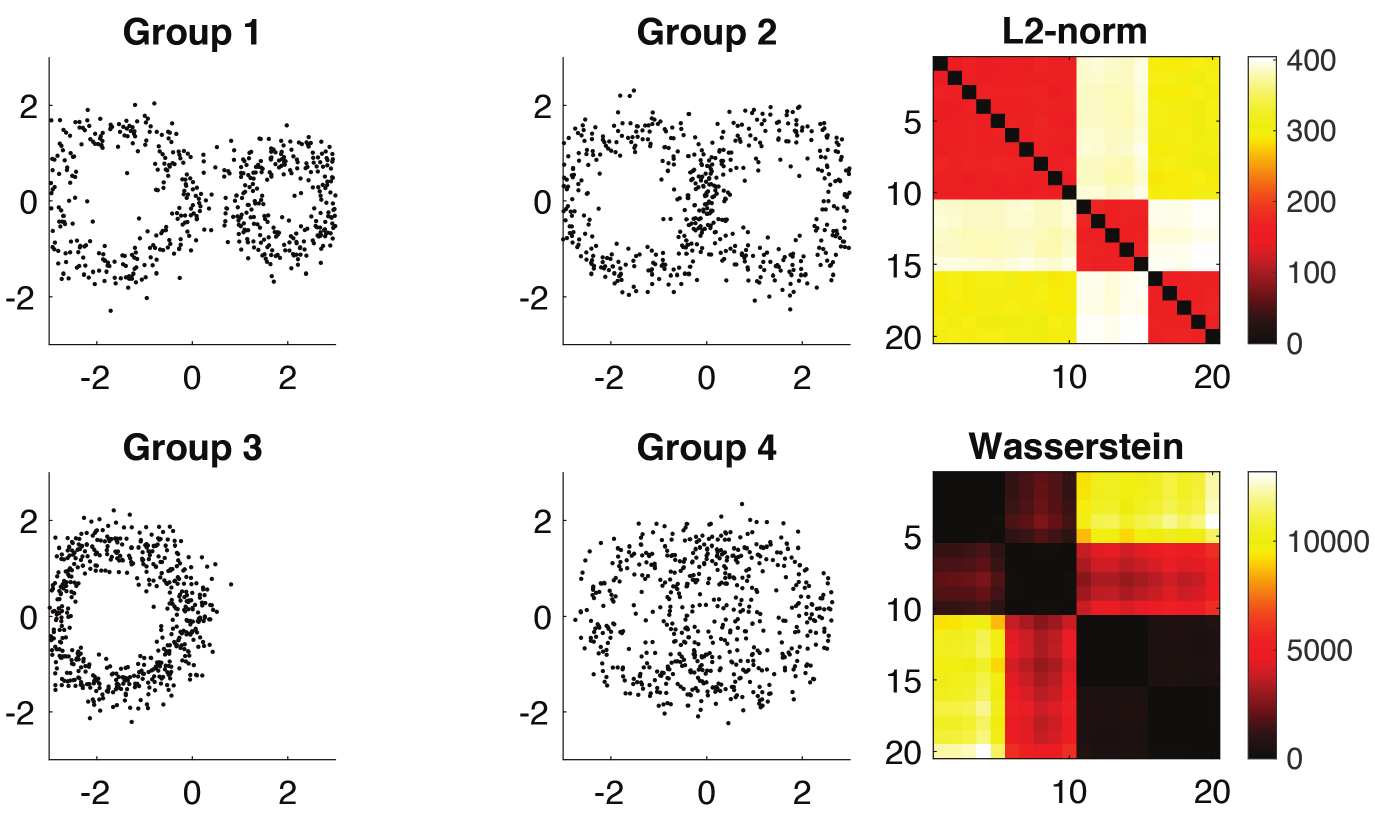}
\caption{Simulation study on  topological difference. The correct clustering method should be able to cluster them since they are all topologically different. Right: the pairwise Euclidean distance ($L_2$-norm) is used in $k$-means and hierarchical clustering. The Wasserstein distance is used in topological clustering.}
\label{fig:simulationdiff}
\end{figure}

We also tested for false negatives when there is topology difference in Figure \ref{fig:simulationdiff}, where all the groups have different number of cycles. All the groups are topologically different and thus we should detect topological differences. The $k$-means clustering achieved $0.83 \pm 0.16$ accuracy. The hierarchical clustering is reporting perfect 1.00 accuracy. On the other hand, the topological clustering achieved respectable $0.98 \pm 0.09$ accuracy.  
It is extremely difficult to separate purely topological signals from geometric signals. Thus, when there is topological difference, it is expected to have geometric signal. Thus, all the methods are expected to perform well. 

Existing clustering methods based on geometric distances will likely to produce significant amount of false positives and and not suitable for topological learning tasks. On the other hand, the proposed Wasserstein distance performed extremely well in both cases and not likely to report false positives or false negatives. The clusterings are performed uing

\begin{verbatim}
acc_WS  = WS_cluster(G)
acc_K   = kmeans_cluster(G)
acc_H   = hierarchical_cluster(G)
\end{verbatim}
\end{example}

\section{Hodge Theory}

Hodge theory provides a unified algebraic framework for analyzing signals defined on simplicial complexes. While classical graph theory studies functions defined on nodes or edges of a graph, Hodge theory extends these ideas to higher-dimensional simplices, including triangles, tetrahedra, and higher-order simplices. This extension enables the decomposition of higher-order signals into orthogonal components associated with global potential, local circulation, and topological cycles \cite{hodge.1989,lim.2020,schaub.2020}.

\subsection{Hodge Laplacian} 

The boundary matrix $\partial_1$ relates nodes to edges is commonly referred as an incidence matrix in graph theory. 
The boundary operator $\partial_k$ can be represented and interpreted as  as the higher dimensional version of {\em incidence matrix} \cite{lee.2018.ISBI,lee.2019.CTIC,schaub.2018}. The standard graph Laplacian can be computed using an incidence matrix as
$$\Delta_0 = \partial_1 \partial_1^{\top},$$
which is also called the 0-th Hodge Laplacian \cite{lee.2018.ISBI}. In general, the $k$-th Hodge Laplacian is defined as
$$\Delta_k =  \partial_{k+1} \partial_{k+1}^{\top} + \partial_k^{\top} \partial_k.$$

The $k$-th Laplacian is a sparse $n_k \times n_k$ positive semi-definite symmetric matrix, where $n_k$ is the number of $k$-simplices   \cite{friedman.1998}. Then the $k$-th Betti number $\beta_k$ is the dimension of $ker \Delta_k$, which is given by computing the rank of $\Delta_k$. The 0th Betti number, the number of connected component, is computed from $\Delta_0$ while the 1st Betti number, the number of independent cycles, is computed from $\Delta_1$. In case of graphs, which is  1-skeletons consisting of only $0$-simplices and $1$-simplies, the boundary matrix $\partial_{2}=0$, thus  the second term in the Hodge Laplacian $\Delta_1$ vanishes and we have  
 \bqn \Delta_1 = \partial_1^{\top}\partial_1. \label{eq:HodgeL1} \eqn
 In brain network studies, brain networks are usually represented as graphs and thus (\ref{eq:HodgeL1}) is more than sufficient unless we model higher order brain connectivity \cite{anand.2023.TMI}. The Hodge Laplacian can be computed differently using the adjacency matrices as
$$\Delta_k = D - A^+ + (k+1)I_{n_k} + A^-,$$
where  $A^+$ and $A^-$ are the upper and lower adjacency matrices between the $k$-simplices. $D=diag(deg(\sigma_1), \cdots, deg(\sigma_{n_k}))$ is the diagonal matrix consisting of the sum of node degrees of simplices $\sigma_j$ \cite{muhammad.2006}. 

The boundary matrices {\tt B} as input, {\tt PH\_hodge.m} outputs the Hodge Laplacian as a cell array {\tt H}. Then {\tt PH\_hodge\_betti.m} computes the Betti numbers through the rank of kernel space of the Hodge Laplacians:
\begin{verbatim}
H=PH_hodge(B)
betti= PH_hodge_betti(H)
\end{verbatim}

\subsection{Hodge decomposition}

\begin{figure}[t]
	\centering
	\includegraphics[width=0.7\linewidth]{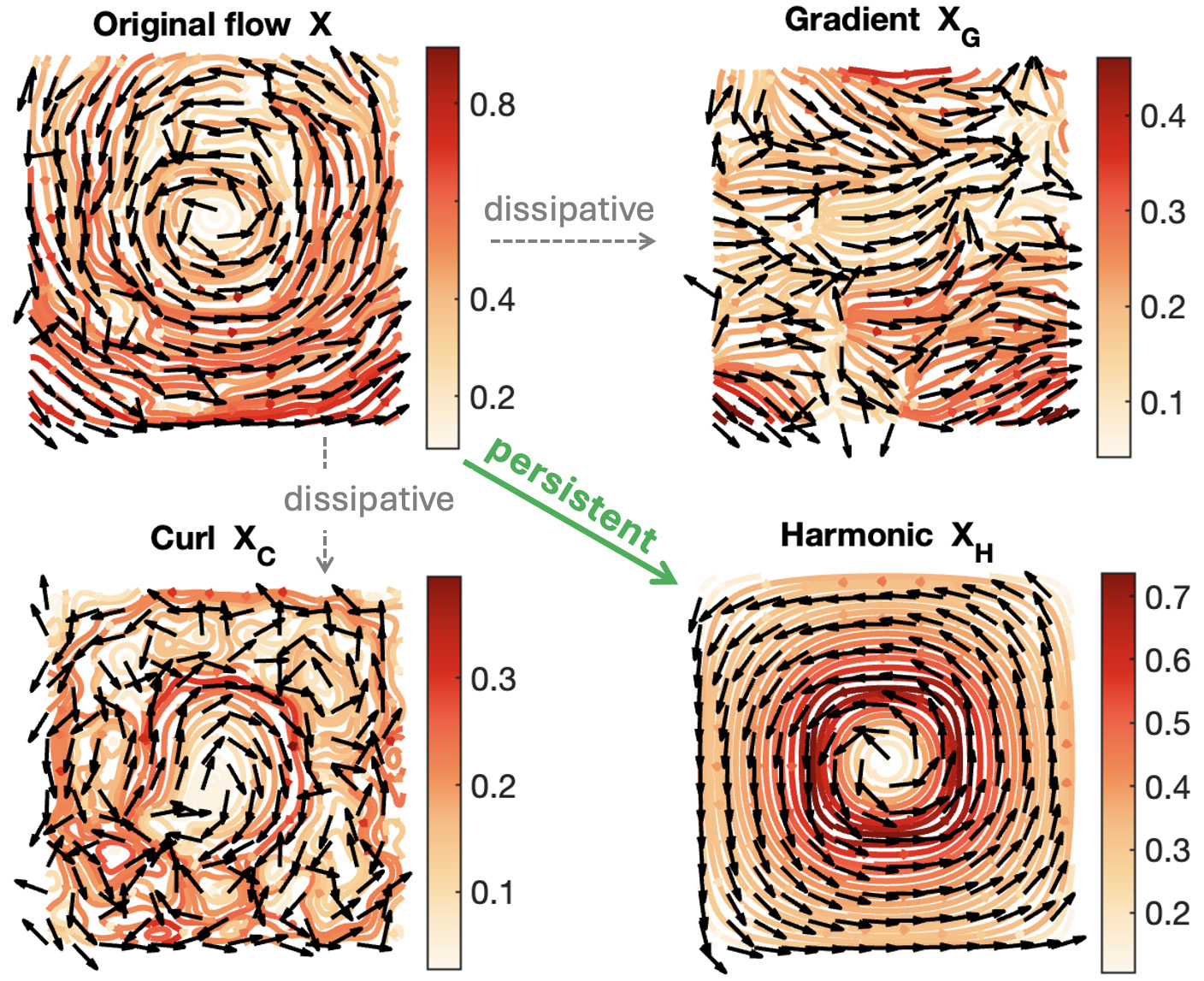}	
\caption{Hodge decomposition of edge signal $X$. The gradient and curl flows $X_G$ and $X_C$ correspond to dissipative structures of the original flow $X$ that decay under Hodge diffusion. In contrast, the harmonic flow $X_H$ is non-dissipative and therefore persists under diffusion, encoding the globally consistent, topologically constrained circulation that remains after transient components dissipate.}
\label{fig:decom}
\end{figure}

Hodge decomposition provides an algebraic framework for decomposing an edge flow into orthogonal components \cite{anand.2024.ISBI,chung.2026.EMBC.causality,chung.2026.ISBI}. Let \(Y\in\mathbb{R}^{p\times p}\) be a directed weighted adjacency matrix represented as an antisymmetric matrix, i.e., \(Y_{ij}=-Y_{ji}\). The corresponding edge flow \(X\in\mathbb{R}^{|E|}\) is obtained by vectorizing the upper triangular entries of \(Y\) according to the edge ordering induced by the simplicial complex. The Hodge decomposition expresses \(X\) as the orthogonal sum
\[
X = X_g + X_c + X_h,
\]
where \(X_g\), \(X_c\), and \(X_h\) denote the gradient, curl, and harmonic components, respectively (Figure \ref{fig:decom}). The gradient component captures source--sink structure induced by node potentials, the curl component represents local circulation around filled triangles, and the harmonic component captures global cyclic flow that cannot be explained by either source--sink structure or local triangular circulation.

The gradient component is the orthogonal projection of \(X\) onto the column space of \(B_1^\top\),
\[
X_g = B_1^\top s,
\]
where the node potential \(s\) is estimated by
\[
\widehat{s}
=
\arg\min_s
\|X-B_1^\top s\|^2.
\]
Similarly, the curl component is the orthogonal projection of \(X\) onto the column space of \(B_2\),
\[
X_c = B_2 z,
\]
where the face potential \(z\) is obtained from
\[
\widehat{z}
=
\arg\min_z
\|X-B_2 z\|^2.
\]
The harmonic component is orthogonal to both the gradient and curl subspaces and satisfies
\[
B_1X_h=0,
\qquad
B_2^\top X_h=0.
\]
Equivalently, the harmonic component is estimated as residual
\[
X_h
=
X-X_g-X_c.
\]
Since
\[
X_h\in\ker(\Delta_1),
\]
the harmonic component represents nontrivial topological cycles that are neither gradients nor boundaries of filled triangles \cite{chung.2026.EMBC.causality}.

\begin{figure}[t]
	\centering
	\includegraphics[width=0.7\linewidth]{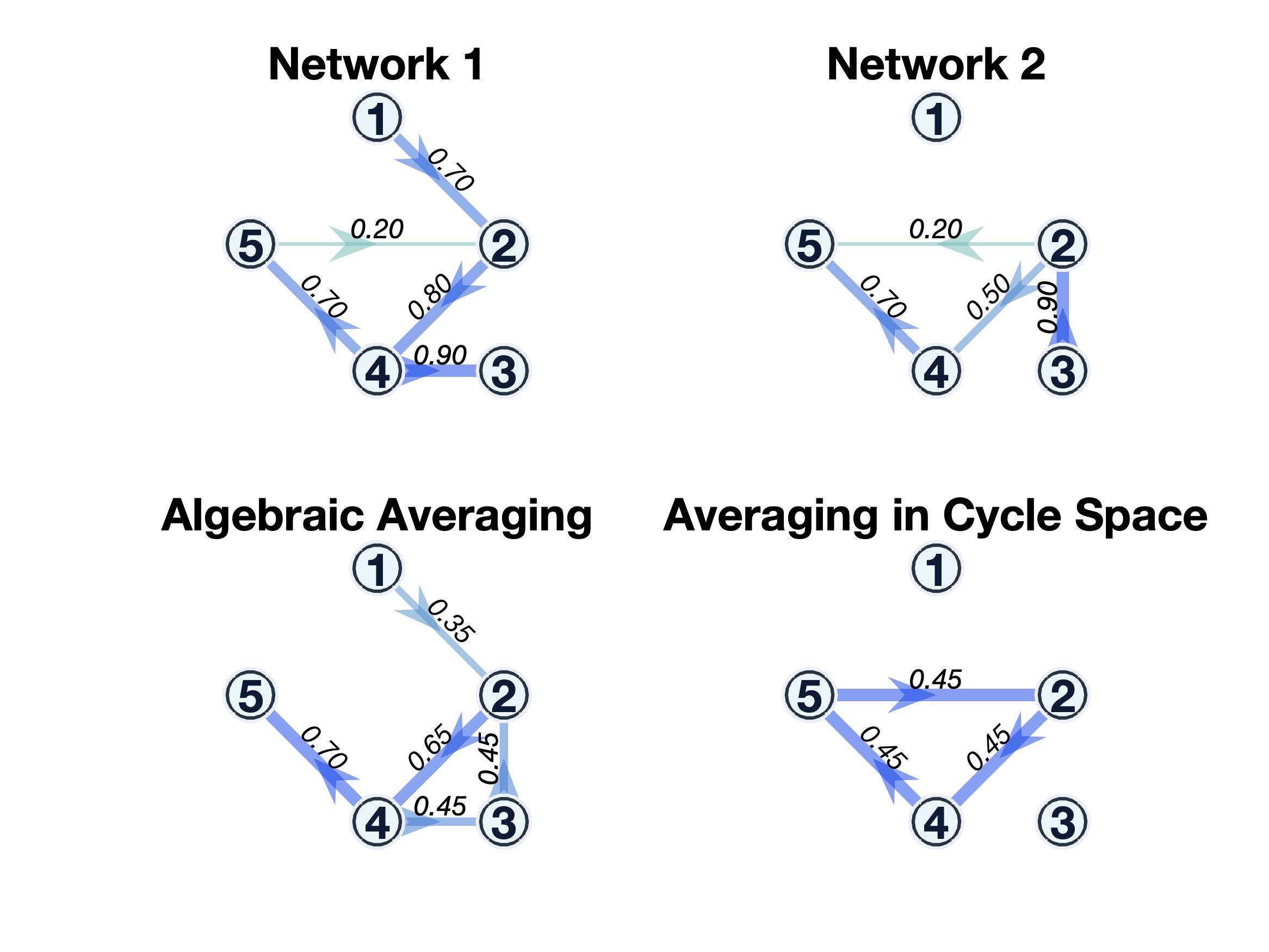}	
\caption{Two networks contain partially inconsistent cyclic structure involving \(2\to4\to5\). Direct algebraic averaging removes this cycle and instead introduces a spurious cycle \(2\to4\to3 \to 2\) that is not present in either network. In contrast, averaging in the harmonic subspace enhances the coherent recurrent cycle \(5\to2\to4\to5\).}
\label{fig:addition}
\end{figure}

The harmonic subspace \(\ker(\Delta_1)\) is a linear subspace of the edge-flow space. If \(X_1,X_2\in\ker(\Delta_1)\) are harmonic flows and \(a,b\in\mathbb{R}\), then
\[
\Delta_1(aX_1+bX_2)
=
a\Delta_1X_1+b\Delta_1X_2
=
0.
\]
Thus, any linear combination of harmonic flows is again harmonic. In particular, adding two admissible cyclic flows produces another admissible cyclic flow. Thus, averaging of cycles are well-defined operations that preserve cyclic consistency (Figure \ref{fig:addition}).

In PH-STAT, the function {\tt Hodge\_vec()} vectorizes the nonzero upper-triangular entries of \(Y\) into an edge-flow vector \(X\) using a lexicographic ordering of node pairs \((i,j)\) with \(i<j\). Since PH-STAT uses sparse matrix representations, zero-weight edges are excluded from both the simplicial complex and the edge flow. Consequently, consistent edge indexing and simplicial support must be maintained across subjects or time points for group-level analyses.

The Hodge decomposition requires the boundary matrices \(B_1\) and \(B_2\), which encode the topology of the underlying simplicial complex. If the data are represented only as a graph, the 1-skeleton is constructed using {\tt Hodge\_1Skeleton()}, and the decomposition contains only the gradient and harmonic components. To compute the curl component, the graph must be extended to a 2-skeleton by specifying which triangular loops are filled as 2-simplices. This is accomplished using {\tt Hodge\_2Skeleton()}, followed by {\tt PH\_boundary()} to construct the boundary matrices. The matrix \(B_2\) encodes the incidence between filled triangles and their boundary edges, so the choice of triangle-filling rule directly determines the curl component.

The decomposition is then computed using
\begin{verbatim}
[Xg, Xc, Xh] = Hodge_decompose(X, B);
\end{verbatim}
The gradient component {\tt Xg} captures source--sink structure induced by node potentials and corresponds to a directed acyclic graph (DAG)-like organization. The curl component {\tt Xc} represents local circulation around filled triangles, while the harmonic component {\tt Xh} captures globally consistent cyclic flow that cannot be explained by either source--sink structure or local triangular circulation.

Each component can be converted back into a MATLAB directed graph format using {\tt edgeflow2digraph()} and visualized with {\tt digraph\_plot()}:

\begin{verbatim}
XG = edgeflow2digraph(kSkeleton, Xg);
XC = edgeflow2digraph(kSkeleton, Xc);
XH = edgeflow2digraph(kSkeleton, Xh);
coord = PH_rips_display([], kSkeleton);
figure; subplot(1,3,1); digraph_plot(XG, coord);
subplot(1,3,2); digraph_plot(XC, coord);
subplot(1,3,3); digraph_plot(XH, coord);
\end{verbatim}

\subsection{Hodge diffusion}

Hodge diffusion generalizes classical graph diffusion from node-valued signals to edge flows \cite{chung.2026.EMBC.causality}. Let
\[
\Delta_1=B_1^\top B_1+B_2B_2^\top
\]
denote the 1-Hodge Laplacian. The evolution of an edge flow \(X(t)\) is governed by the diffusion equation
\[
\frac{dX(t)}{dt}
=
-\Delta_1X(t),
\qquad
X(0)=X.
\]
The solution is
$
X(t)
=
e^{-t\Delta_1}X,
$
which continuously smooths the edge flow while respecting the topology of the underlying simplicial complex.

The diffusion process admits a simple spectral representation. Let
\[
\Delta_1\phi_l
=
\lambda_l\phi_l,
\]
where \(\phi_l\) is an eigenvector and \(\lambda_l\ge0\) is the corresponding eigenvalue. Then
\[
X(t)
=
\sum_l
e^{-t\lambda_l}
c_l
\phi_l,
\qquad
c_l=\phi_l^\top X.
\]
Since eigenmodes associated with larger eigenvalues decay more rapidly, gradient and curl components gradually disappear during diffusion. In contrast, the eigenvectors associated with zero eigenvalues span the harmonic subspace $\ker(\Delta_1),$
whose dimension equals the first Betti number \(\beta_1\). Consequently,
\[
\lim_{t\rightarrow\infty}
X(t)
=
X_h,
\]
where \(X_h\) is the harmonic component of the Hodge decomposition (Figure \ref{fig:decom}). Thus, Hodge diffusion provides an alternative characterization of harmonic flow as the steady-state solution of edge-flow diffusion. If $\Phi_0=
[\phi_1,\ldots,\phi_{\beta_1}]$
contains the eigenvectors corresponding to zero eigenvalues, then the harmonic component is the orthogonal projection
\[
X_h
=
\Phi_0\Phi_0^\top X.
\]

In PH-STAT, the 1-Hodge Laplacian and its eigensystem are computed as
\begin{verbatim}
Delta1 = Hodge_laplacian(B);
[Phi,Lambda] = eig(full(Delta1));
lambda = diag(Lambda);
idx0 = find(abs(lambda) < 1e-8);
Phi0 = Phi(:,idx0);
\end{verbatim}
The harmonic component can then be obtained directly from the spectral projection
\begin{verbatim}
Xh = Phi0 * (Phi0' * X);
\end{verbatim}
or equivalently from the Hodge decomposition
\begin{verbatim}
[Xg,Xc,Xh] = Hodge_decompose(X,B);
\end{verbatim}
The edge flow after diffusion time \(t\) is computed by
\begin{verbatim}
t = 1;
c  = Phi' * X;
Xt = Phi * (exp(-t*lambda) .* c);
\end{verbatim}

\section*{Acknowledgement}
The project is supported by NIH R01 EB028753 and  NSF MDS-2010778. We also like to thank Hyekyung Lee of Seoul National University and Tananun Songdechakraiwut of University of Wisconsin-Madison, Vijay D. Anand of University College of London for the contribution of some of Matlab functions in PH-STAT package.


\bibliographystyle{splncs04}

\begin{thebibliography}{10}
\providecommand{\url}[1]{\texttt{#1}}
\providecommand{\urlprefix}{URL }
\providecommand{\doi}[1]{https://doi.org/#1}

\bibitem{adler.2010}
Adler, R., Bobrowski, O., Borman, M., Subag, E., Weinberger, S.: Persistent
  homology for random fields and complexes. In: Borrowing strength: theory
  powering applications--a Festschrift for Lawrence D. Brown, pp. 124--143.
  Institute of Mathematical Statistics (2010)

\bibitem{agueh.2011}
Agueh, M., Carlier, G.: Barycenters in the {W}asserstein space. SIAM Journal on
  Mathematical Analysis  \textbf{43},  904--924 (2011)

\bibitem{anand.2023.TMI}
Anand, D., Chung, M.: {Hodge-Laplacian} of brain networks. IEEE Transactions on
  Medical Imaging  \textbf{42},  1563--1573 (2023)

\bibitem{anand.2024.ISBI}
Anand, D., Chung, M.: Hodge-decomposition of brain networks. In: 2024 IEEE
  International Symposium on Biomedical Imaging (ISBI). pp.~1--5. IEEE (2024)

\bibitem{anirudh.2016}
Anirudh, R., Thiagarajan, J., Kim, I., Polonik, W.: Autism spectrum disorder
  classification using graph kernels on multidimensional time series. arXiv
  preprint arXiv:1611.09897  (2016)

\bibitem{bendich.2016}
Bendich, P., Marron, J., Miller, E., Pieloch, A., Skwerer, S.: Persistent
  homology analysis of brain artery trees. The annals of applied statistics
  \textbf{10}, ~198 (2016)

\bibitem{berwald.2018}
Berwald, J., Gottlieb, J., Munch, E.: Computing wasserstein distance for
  persistence diagrams on a quantum computer. arXiv:1809.06433  (2018)

\bibitem{bishop.2006}
Bishop, C.: Pattern recognition and machine learning. Springer (2006)

\bibitem{cai.2020}
Cai, Y., Zhang, J., Xiao, T., Peng, H., Sterling, S., Walsh, R., Rawson, S.,
  Rits-Volloch, S., Chen, B.: Distinct conformational states of {SARS-CoV-2}
  spike protein. Science  \textbf{369},  1586--1592 (2020)

\bibitem{canas.2012}
Canas, G., Rosasco, L.: Learning probability measures with respect to optimal
  transport metrics. arXiv preprint arXiv:1209.1077  (2012)

\bibitem{carlsson.2008}
Carlsson, G., Memoli, F.: Persistent clustering and a theorem of {J.
  Kleinberg}. arXiv preprint arXiv:0808.2241  (2008)

\bibitem{cassidy.2015}
Cassidy, B., Rae, C., Solo, V.: Brain activity: conditional dissimilarity and
  persistent homology. In: IEEE 12th International Symposium on Biomedical
  Imaging (ISBI). pp. 1356--1359 (2015)

\bibitem{chung.2011.SPIE}
Chung, M., Adluru, N., Dalton, K., Alexander, A., Davidson, R.: Scalable brain
  network construction on white matter fibers. In: Proc. of SPIE. vol.~7962, p.
  79624G (2011)

\bibitem{chung.2026.EMBC.causality}
Chung, M., Anand, D., El-Yaagoubi, A., Jung, J.H., Qiu, A., Ombao, H.:
  Causality as a minimum energy principle. In: Proceedings of the Annual
  International Conference of the IEEE Engineering in Medicine and Biology
  Society (EMBC). pp.~1--5 (2026)

\bibitem{chung.2009.IPMI}
Chung, M., Bubenik, P., Kim, P.: Persistence diagrams of cortical surface data.
  Proceedings of the 21st International Conference on Information Processing in
  Medical Imaging (IPMI), Lecture Notes in Computer Science (LNCS)
  \textbf{5636},  386--397 (2009)

\bibitem{chung.2007.TMI}
Chung, M., Dalton, K., Shen, L., Evans, A., Davidson, R.: Weighted {Fourier}
  representation and its application to quantifying the amount of gray matter.
  IEEE Transactions on Medical Imaging  \textbf{26},  566--581 (2007)

\bibitem{chung.2013.MICCAI}
Chung, M., Hanson, J., Lee, H., Adluru, N., Alexander, A.L., Davidson, R.,
  Pollak, S.: Persistent homological sparse network approach to detecting white
  matter abnormality in maltreated children: {MRI} and {DTI} multimodal study.
  MICCAI, Lecture Notes in Computer Science (LNCS)  \textbf{8149},  300--307
  (2013)

\bibitem{chung.2015.TMI}
Chung, M., Hanson, J., Ye, J., Davidson, R., Pollak, S.: Persistent homology in
  sparse regression and its application to brain morphometry. IEEE Transactions
  on Medical Imaging  \textbf{34},  1928--1939 (2015)

\bibitem{chung.2019.ISBI}
Chung, M., Huang, S.G., Gritsenko, A., Shen, L., Lee, H.: Statistical inference
  on the number of cycles in brain networks. In: 2019 IEEE 16th International
  Symposium on Biomedical Imaging (ISBI 2019). pp. 113--116. IEEE (2019)

\bibitem{chung.2019.NN}
Chung, M., Lee, H., DiChristofano, A., Ombao, H., Solo, V.: Exact topological
  inference of the resting-state brain networks in twins. Network Neuroscience
  \textbf{3},  674--694 (2019)

\bibitem{chung.2026.ISBI}
Chung, M., Maccotta, L., Struck, A.: Counterfactual analysis of brain network
  dynamics. In: Proceedings of the IEEE International Symposium on Biomedical
  Imaging (ISBI). pp.~1--5 (2026)

\bibitem{chung.2021.MICCAI}
Chung, M., Ombao, H.: Lattice paths for persistent diagrams. In:
  Interpretability of Machine Intelligence in Medical Image Computing, and
  Topological Data Analysis and Its Applications for Medical Data, LNCS 12929,
  pp. 77--86 (2021)

\bibitem{chung.2023.NI}
Chung, M., Ramos, C., De~Paiva, F., Mathis, J., Prabhakaran, V., Nair, V.,
  Meyerand, M., Hermann, B., Binder, J., Struck, A.: Unified topological
  inference for brain networks in temporal lobe epilepsy using the
  {Wasserstein} distance. NeuroImage  \textbf{284},  120436 (2023)

\bibitem{chung.2020.arXiv}
Chung, M., Smith, A., Shiu, G.: Reviews: Topological distances and losses for
  brain networks. arXiv e-prints pp. arXiv--2102.08623 (2020),
  \url{https://arxiv.org/pdf/2102.08623}

\bibitem{chung.2017.IPMI}
Chung, M., Vilalta-Gil, V., Lee, H., Rathouz, P., Lahey, B., Zald, D.: Exact
  topological inference for paired brain networks via persistent homology.
  Information Processing in Medical Imaging (IPMI), Lecture Notes in Computer
  Science (LNCS)  \textbf{10265},  299--310 (2017)

\bibitem{chung.2019.CNI}
Chung, M., Xie, L., Huang, S.G., Wang, Y., Yan, J., Shen, L.: Rapid
  acceleration of the permutation test via transpositions. International
  Workshop on Connectomics in Neuroimaging  \textbf{11848},  42--53 (2019)

\bibitem{clough.2019}
Clough, J., Oksuz, I., Byrne, N., Schnabel, J., King, A.: Explicit topological
  priors for deep-learning based image segmentation using persistent homology.
  In: International Conference on Information Processing in Medical Imaging.
  pp. 16--28. Springer (2019)

\bibitem{cohensteiner.2007}
Cohen-Steiner, D., Edelsbrunner, H., Harer, J.: Stability of persistence
  diagrams. Discrete and Computational Geometry  \textbf{37},  103--120 (2007)

\bibitem{cuturi.2014}
Cuturi, M., Doucet, A.: Fast computation of {W}asserstein barycenters. In:
  International conference on machine learning. pp. 685--693. PMLR (2014)

\bibitem{dubey.2019}
Dubey, P., M{\"u}ller, H.G.: Fr{\'e}chet analysis of variance for random
  objects. Biometrika  \textbf{106},  803--821 (2019)

\bibitem{duff.2001}
Duff, I., Koster, J.: On algorithms for permuting large entries to the diagonal
  of a sparse matrix. SIAM Journal on Matrix Analysis and Applications
  \textbf{22}(4),  973--996 (2001)

\bibitem{edelsbrunner.2008}
Edelsbrunner, H., Harer, J.: Persistent homology - a survey. Contemporary
  Mathematics  \textbf{453},  257--282 (2008)

\bibitem{edelsbrunner.2010}
Edelsbrunner, H., Harer, J.: Computational topology: {A}n introduction.
  American Mathematical Society (2010)

\bibitem{edelsbrunner.2002}
Edelsbrunner, H., Letscher, D., Zomorodian, A.: Topological persistence and
  simplification. Discrete and Computational Geometry  \textbf{28},  511--533
  (2002)

\bibitem{edmonds.1972}
Edmonds, J., Karp, R.: Theoretical improvements in algorithmic efficiency for
  network flow problems. Journal of the ACM (JACM)  \textbf{19},  248--264
  (1972)

\bibitem{feller.2008}
Feller, W.: An introduction to probability theory and its applications, vol.~2.
  John Wiley \& Sons (2008)

\bibitem{friedman.1998}
Friedman, J.: Computing betti numbers via combinatorial laplacians.
  Algorithmica  \textbf{21}(4),  331--346 (1998)

\bibitem{friston.2002}
Friston., K.: A short history of statistical parametric mapping in functional
  neuroimaging. Tech. Rep. Technical report, Wellcome Department of Imaging
  Neuroscience, ION, UCL., London, UK. (2002)

\bibitem{ghrist.2008}
Ghrist, R.: Barcodes: The persistent topology of data. Bulletin of the American
  Mathematical Society  \textbf{45},  61--75 (2008)

\bibitem{guo.2020}
Guo, X., Srivastava, A.: Representations, metrics and statistics for shape
  analysis of elastic graphs. In: Proceedings of the IEEE/CVF Conference on
  Computer Vision and Pattern Recognition Workshops. pp. 832--833 (2020)

\bibitem{hart.1999}
Hart, J.: Computational topology for shape modeling. In: Proceedings of the
  International Conference on Shape Modeling and Applications. pp. 36--43
  (1999)

\bibitem{hatcher.2002}
Hatcher, A.: Algebraic topology. Cambridge University Press (2002)

\bibitem{hodge.1989}
Hodge, W.: The theory and applications of harmonic integrals. CUP Archive
  (1989)

\bibitem{hu.2019}
Hu, X., Li, F., Samaras, D., Chen, C.: Topology-preserving deep image
  segmentation. In: Advances in Neural Information Processing Systems. pp.
  5657--5668 (2019)

\bibitem{huang.2020.NM}
Huang, S.G., Samdin, S.T., Ting, C., Ombao, H., Chung, M.: Statistical model
  for dynamically-changing correlation matrices with application to brain
  connectivity. Journal of Neuroscience Methods  \textbf{331},  108480 (2020)

\bibitem{khalid.2014}
Khalid, A., Kim, B., Chung, M., Ye, J., Jeon, D.: Tracing the evolution of
  multi-scale functional networks in a mouse model of depression using
  persistent brain network homology. NeuroImage  \textbf{101},  351--363 (2014)

\bibitem{kiebel.1996}
Kiebel, S., Poline, J.P., Friston, K., Holmes, A., Worsley, K.: Robust
  smoothness estimation in statistical parametric maps using standardized
  residuals from the general linear model. NeuroImage  \textbf{10},  756--766
  (1999)

\bibitem{kullback.1951}
Kullback, S., Leibler, R.: On information and sufficiency. The Annals of
  Mathematical Statistics  \textbf{22},  79--86 (1951)

\bibitem{le.2000}
Le, H., Kume, A.: The {Fr{\'e}chet} mean shape and the shape of the means.
  Advances in Applied Probability  \textbf{32},  101--113 (2000)

\bibitem{lee.2014.MICCAI}
Lee, H., Chung, M.K., K.H., Lee, D.: Hole detection in metabolic connectivity
  of {Alzheimer's} disease using k-{Laplacian}. MICCAI, Lecture Notes in
  Computer Science  \textbf{8675},  297--304 (2014)

\bibitem{lee.2019.CTIC}
Lee, H., Chung, M., Choi, H., K., H., Ha, S., Kim, Y., Lee, D.: Harmonic holes
  as the submodules of brain network and network dissimilarity. International
  Workshop on Computational Topology in Image Context, Lecture Notes in
  Computer Science pp. 110--122 (2019)

\bibitem{lee.2018.ISBI}
Lee, H., Chung, M., Kang, H., Choi, H., Kim, Y., Lee, D.: Abnormal hole
  detection in brain connectivity by kernel density of persistence diagram and
  {Hodge Laplacian}. In: {IEEE International Symposium on Biomedical Imaging
  (ISBI)}. pp. 20--23 (2018)

\bibitem{lee.2011.MICCAI}
Lee, H., Chung, M., Kang, H., Kim, B.N., Lee, D.: Computing the shape of brain
  networks using graph filtration and {Gromov-Hausdorff} metric. {MICCAI,
  Lecture Notes in Computer Science}  \textbf{6892},  302--309 (2011)

\bibitem{lee.2011.ISBI}
Lee, H., Chung, M., Kang, H., Kim, B.N., Lee, D.: Discriminative persistent
  homology of brain networks. In: {IEEE International Symposium on Biomedical
  Imaging (ISBI)}. pp. 841--844 (2011)

\bibitem{lee.2012.TMI}
Lee, H., Kang, H., Chung, M., Kim, B.N., Lee, D.: Persistent brain network
  homology from the perspective of dendrogram. IEEE Transactions on Medical
  Imaging  \textbf{31},  2267--2277 (2012)

\bibitem{lee.2017.HBM}
Lee, H., Kang, H., Chung, M., Lim, S., Kim, B.N., Lee, D.: Integrated
  multimodal network approach to {PET} and {MRI} based on multidimensional
  persistent homology. Human Brain Mapping  \textbf{38},  1387--1402 (2017)

\bibitem{lee.2018.deep}
Lee, M., Xiong, Y., Yu, G., Li, G.Y.: Deep neural networks for linear sum
  assignment problems. IEEE Wireless Communications Letters  \textbf{7},
  962--965 (2018)

\bibitem{li.2017}
Li, Y., Wang, D., Ascoli, G., Mitra, P., Wang, Y.: Metrics for comparing
  neuronal tree shapes based on persistent homology. PloS one  \textbf{12}(8),
  e0182184 (2017)

\bibitem{lim.2020}
Lim, L.H.: Hodge laplacians on graphs. Siam Review  \textbf{62},  685--715
  (2020)

\bibitem{mi.2018}
Mi, L., Zhang, W., Gu, X., Wang, Y.: Variational wasserstein clustering. In:
  Proceedings of the European Conference on Computer Vision (ECCV). pp.
  322--337 (2018)

\bibitem{miller.1997}
Miller, M., Banerjee, A., Christensen, G., Joshi, S., Khaneja, N., Grenander,
  U., Matejic, L.: Statistical methods in computational anatomy. Statistical
  Methods in Medical Research  \textbf{6},  267--299 (1997)

\bibitem{milnor.1973}
Milnor, J.: Morse Theory. Princeton University Press (1973)

\bibitem{morozov.2008}
Morozov, D.: Homological Illusions of Persistence and Stability. Ph.D. thesis,
  Duke University (2008)

\bibitem{muhammad.2006}
Muhammad, A., Egerstedt, M.: Control using higher order laplacians in network
  topologies. In: Proc. of 17th International Symposium on Mathematical Theory
  of Networks and Systems. pp. 1024--1038 (2006)

\bibitem{naiman.1990}
Naiman, D.: volumes for tubular neighborhoods of spherical polyhedra and
  statistical inference. Ann. Statist.  \textbf{18},  685--716 (1990)

\bibitem{nichols.2003}
Nichols, T., Hayasaka, S.: Controlling the familywise error rate in functional
  neuroimaging: a comparative review. Stat Methods Med. Res.  \textbf{12},
  419--446 (2003)

\bibitem{otter.2017}
Otter, N., Porter, M., Tillmann, U., Grindrod, P., Harrington, H.: A roadmap
  for the computation of persistent homology. EPJ Data Science  \textbf{6}(1),
  ~17 (2017)

\bibitem{palande.2017}
Palande, S., Jose, V., Zielinski, B., Anderson, J., Fletcher, P., Wang, B.:
  Revisiting abnormalities in brain network architecture underlying autism
  using topology-inspired statistical inference (2017)

\bibitem{petri.2014}
Petri, G., Expert, P., Turkheimer, F., Carhart-Harris, R., Nutt, D., Hellyer,
  P., Vaccarino, F.: Homological scaffolds of brain functional networks.
  Journal of The Royal Society Interface  \textbf{11},  20140873 (2014)

\bibitem{pothen.1990}
Pothen, A., Fan, C.: {Computing the block triangular form of a sparse matrix}.
  ACM Transactions on Mathematical Software (TOMS)  \textbf{16}, ~324 (1990)

\bibitem{rabin.2011}
Rabin, J., Peyr{\'e}, G., Delon, J., Bernot, M.: Wasserstein barycenter and its
  application to texture mixing. In: International Conference on Scale Space
  and Variational Methods in Computer Vision. pp. 435--446. Springer (2011)

\bibitem{schaub.2018}
Schaub, M., Benson, A., Horn, P., Lippner, G., Jadbabaie, A.: Random walks on
  simplicial complexes and the normalized hodge laplacian. arXiv preprint
  arXiv:1807.05044  (2018)

\bibitem{schaub.2020}
Schaub, M., Benson, A., Horn, P., Lippner, G., Jadbabaie, A.: Random walks on
  simplicial complexes and the normalized {Hodge} 1-laplacian. SIAM Review
  \textbf{62},  353--391 (2020)

\bibitem{sheehy.2013}
Sheehy, D.: Linear-size approximations to the {Vietoris--Rips} filtration.
  Discrete \& Computational Geometry  \textbf{49},  778--796 (2013)

\bibitem{deSilva.2007}
de~Silva, V., Ghrist, R.: Homological sensor networks. Notic Amer Math Soc
  \textbf{54},  10--17 (2007)

\bibitem{solo.2018}
Solo, V., Poline, J., Lindquist, M., Simpson, S., Bowman, D., C.M., Cassidy,
  B.: Connectivity in {fMRI}: Blind spots and breakthroughs. IEEE Transactions
  on Medical Imaging  \textbf{37},  1537--1550 (2018)

\bibitem{song.2023}
Songdechakraiwut, T., Chung, M.: Topological learning for brain networks.
  Annals of Applied Statistics  \textbf{17},  403--433 (2023)

\bibitem{song.2021.MICCAI}
Songdechakraiwut, T., Shen, L., Chung, M.: Topological learning and its
  application to multimodal brain network integration. Medical Image Computing
  and Computer Assisted Intervention (MICCAI)  \textbf{12902},  166--176 (2021)

\bibitem{taylor.2008}
Taylor, J., Worsley, K.: Random fields of multivariate test statistics, with
  applications to shape analysis. Annals of Statistics  \textbf{36},  1--27
  (2008)

\bibitem{topaz.2015}
Topaz, C., Ziegelmeier, L., Halverson, T.: Topological data analysis of
  biological aggregation models. PLoS One p. e0126383 (2015)

\bibitem{turner.2014}
Turner, K., Mileyko, Y., Mukherjee, S., Harer, J.: Fr{\'e}chet means for
  distributions of persistence diagrams. Discrete \& Computational Geometry
  \textbf{52},  44--70 (2014)

\bibitem{vallender.1974}
Vallender, S.: Calculation of the {W}asserstein distance between probability
  distributions on the line. Theory of Probability \& Its Applications
  \textbf{18},  784--786 (1974)

\bibitem{walls.2020}
Walls, A., Park, Y.J., Tortorici, M., Wall, A., McGuire, A., Veesler, D.:
  Structure, function, and antigenicity of the {SARS-CoV-2} spike glycoprotein.
  Cell  \textbf{181},  281--292 (2020)

\bibitem{wang.2018.annals}
Wang, Y., Ombao, H., Chung, M.: Topological data analysis of single-trial
  electroencephalographic signals. Annals of Applied Statistics  \textbf{12},
  1506--1534 (2018)

\bibitem{wong.2016}
Wong, E., Palande, S., Wang, B., Zielinski, B., Anderson, J., Fletcher, P.:
  Kernel partial least squares regression for relating functional brain network
  topology to clinical measures of behavior. In: IEEE International Symposium
  on Biomedical Imaging (ISBI). pp. 1303--1306 (2016)

\bibitem{worsley.1996}
Worsley, K., Marrett, S., Neelin, P., Vandal, A., Friston, K., Evans, A.: A
  unified statistical approach for determining significant signals in images of
  cerebral activation. Human Brain Mapping  \textbf{4},  58--73 (1996)

\bibitem{yang.2020}
Yang, Z., Wen, J., Davatzikos, C.: Smile-{GANs}: Semi-supervised clustering via
  {GANs} for dissecting brain disease heterogeneity from medical images. arXiv
  preprint  \textbf{arXiv},  2006.15255 (2020)

\bibitem{zemel.2019}
Zemel, Y., Panaretos, V.: Fr{\'e}chet means and procrustes analysis in
  {Wasserstein} space. Bernoulli  \textbf{25},  932--976 (2019)

\bibitem{zomorodian.2001}
Zomorodian, A.: Computing and Comprehending Topology: Persistence and
  Hierarchical Morse Complexes. Ph.D. Thesis, University of Illinois,
  Urbana-Champaign (2001)

\bibitem{zomorodian.2009}
Zomorodian, A.: Topology for computing. Cambridge University Press, Cambridge
  (2009)

\bibitem{zomorodian.2005}
Zomorodian, A., Carlsson, G.: Computing persistent homology. Discrete and
  Computational Geometry  \textbf{33},  249--274 (2005)

\end{thebibliography}

\end{document}